\documentclass{amsart}
\usepackage{comment}
\usepackage{MnSymbol}
\usepackage{extarrows}
\usepackage{lscape}
\usepackage[all,cmtip]{xy}

\DeclareMathOperator{\Spec}{Spec}
\DeclareMathOperator{\Sym}{Sym}
\newcommand{\pdual}{\check{\mathbb{P}}^2}
\DeclareMathOperator{\PGL}{PGL}
\DeclareMathOperator{\PSL}{PSL}
\DeclareMathOperator{\Aut}{Aut}

\DeclareMathOperator{\GL}{GL}
\DeclareMathOperator{\SL}{SL}

\DeclareMathOperator{\Pic}{Pic}
\DeclareMathOperator{\red}{red}

\DeclareMathOperator{\divisor}{div}

\DeclareMathOperator{\Gal}{Gal}

\newcommand{\Q}{{\mathbb Q}}
\newcommand{\Z}{{\mathbb Z}}
\newcommand{\C}{{\mathbb C}}
\newcommand{\R}{{\mathbb R}}
\newcommand{\F}{{\mathbb F}}
\newcommand{\PP}{{\mathbb P}}
\newcommand{\cA}{\mathcal{A}}

\newcommand{\cM}{\mathcal{M}}
\newcommand{\cN}{\mathcal{N}}

\newcommand{\OO}{{\mathcal O}}

\newcommand{\ff}{\mathfrak{f}}

\newcommand{\fp}{\mathfrak{p}}
\newcommand{\fq}{\mathfrak{q}}

\begin {document}

\newtheorem{thm}{Theorem}
\newtheorem{lem}{Lemma}[section]
\newtheorem{prop}[lem]{Proposition}

\newtheorem{cor}[lem]{Corollary}

\theoremstyle{definition}

\theoremstyle{remark}

\title[]{Elliptic Curves over Real Quadratic Fields\\
are Modular}

\author{Nuno Freitas, Bao V. Le Hung, Samir Siksek}
\address{
 Mathematisches Institut\\
Universit\"{a}t Bayreuth\\
95440 Bayreuth, Germany
}
\email{nunobfreitas@gmail.com}

\address{
Department of Mathematics\\ 
Harvard University\\ 
One Oxford Street\\ 
Cambridge\\ MA 02138, USA}

\email{lhvietbao@googlemail.com}

\address{Mathematics Institute\\
	University of Warwick\\
	CV4 7AL \\
	United Kingdom}

\email{samir.siksek@gmail.com}

\date{\today}
\thanks{
The first-named author is supported
through a grant within the framework of the DFG Priority Programme 1489
{\em Algorithmic and Experimental Methods in Algebra, Geometry and Number Theory}.
The third-named 
author is supported by an EPSRC Leadership Fellowship EP/G007268/1,
and EPSRC {\em LMF: L-Functions and Modular Forms} Programme Grant
EP/K034383/1.
}


\begin{abstract}
We prove that all elliptic curves defined over real quadratic fields are modular.
\end{abstract}
\maketitle

\makeatletter
\renewcommand{\tocsection}[3]{%
  \indentlabel{\@ifnotempty{#2}{\ignorespaces#1 #2.\quad}}#3\dotfill
}

\makeatother
\tableofcontents

\part{Introduction}

\section{Summary of Results}

One of the great achievements of modern number theory is the proof 
by
Breuil, Conrad, Diamond and Taylor ~\cite{modularity}
of the
Shimura--Taniyama (or modularity) conjecture for elliptic curves over $\Q$.
This built
on the earlier
pioneering work of Wiles ~\cite{Wiles}, and of Taylor and Wiles ~\cite{TW}, 
which established the modularity
of semistable elliptic curves over $\Q$, and was historically the final
step in the proof of  
Fermat's Last Theorem. In this paper we prove the analogous modularity
statement for elliptic curves over real quadratic fields.

\bigskip

Let $K$ be a totally real number field, and write $G_K:=\Gal(\overline{\Q}/K)$
for its absolute Galois group. Let $E$ be an elliptic curve over $K$.
Recall that $E$ 
is \textbf{modular} if 
there exists a Hilbert cuspidal eigenform $\ff$
over $K$ of parallel weight $2$, with rational
Hecke eigenvalues, such that the Hasse--Weil
$\mathrm{L}$-function of $E$ is equal to the $\mathrm{L}$-function
of $\ff$. A more conceptual way to phrase this is that there is an isomorphism of compatible systems of Galois representations 
\[ 
\rho_{E,p} \cong \rho_{\ff,p}  
\]
where the left-hand side is the Galois representation arising from the action
of $G_K$ on the $p$-adic Tate module $T_p(E)$, while the right-hand side is the
Galois representation associated to $\ff$ by Eichler and Shimura (when $K=\Q$)
~\cite{DS}, and by Carayol ~\cite{Carayol1} ~\cite{Carayol2}, Blasius and Rogawski
~\cite{Blasius-Rogawski}, Wiles~\cite{Wiles2}
and Taylor ~\cite{TaylorHilbert} (for general totally
real $K$). 

Work of Jarvis and Manoharmayum \cite{JM} establishes modularity of semistable elliptic 
curves over $\Q(\sqrt{2})$ and $\Q(\sqrt{17})$, and there are various other works 
\cite{Allen},
\cite{DF2}, 
which establish modularity under local assumptions on the curve $E$ and the field $K$.
In this paper, we prove modularity of all elliptic curves over all real quadratic fields.
\begin{thm}\label{thm:modquadratic}
Let $E$ be an elliptic curve over a real quadratic field $K$.
Then
$E$ is modular.
\end{thm}

Before sketching the strategy of the proof of Theorem~\ref{thm:modquadratic},
we need a little more notation and terminology.
Let $p$ be a rational prime. 
Denote by 
\[
\overline{\rho}_{E,p} \; : \; G_K \rightarrow \Aut(E[p]) \cong \GL_2(\F_p)
\] 
the representation
giving the action of $G_K$ on the $p$-torsion of $E$.
We say $\overline{\rho}_{E,p}$ is {\bf modular} if
there exists a Hilbert cuspidal eigenform $\ff$
over $K$ of parallel weight $2$, 
and a place $\lambda \mid p$ of $\overline{\Q}$
such that $\overline{\rho}^{\mathrm{ss}}_{E,p} \sim \overline{\rho}^{\mathrm{ss}}_{\ff,\lambda}$.  

\bigskip

For our purposes, the following is a convenient sketch of the main steps in the proof  
(\cite{Wiles}, \cite{TW},
\cite{Diamond}, \cite{CDT}, \cite{modularity}) 
of modularity of elliptic curves
over $\Q$:
\begin{enumerate}
  \item[Step I] (modularity lifting): one shows that 
if
$\overline{\rho}_{E,p}$ is modular and the image of $\overline{\rho}_{E,p}$ is
\lq sufficiently large\rq\ (in a suitable sense), then $\rho_{E,p}$ is modular.
  \item[Step II] (Langlands--Tunnell): if $\overline{\rho}_{E,3}$ is 
irreducible then it is modular.
 \item[Step III] (3--5 modularity switching): with the help of Steps I and II, one
shows that if $\overline{\rho}_{E,5}$ is irreducible then it is modular.
  \item[Step IV] (dealing with exceptions): putting Steps I, II, III together,
one sees that $E$ is modular except possibly if the image of $\overline{\rho}_{E,p}$
is \lq small\rq\ for $p=3$ and $5$ simultaneously. Such elliptic curves
have certain mod $3$ and mod $5$ level structures, and so give
rise to rational points on certain modular curves. 
By 
enumerating these rational points (of which there are only a handful)
one shows that these elliptic curves are modular by some other means.
\end{enumerate}
We are ready to outline the proof of Theorem~\ref{thm:modquadratic}
and its similarities and differences from the above;
we also outline the contents of the paper. 
First we need a suitable modularity lifting theorem (Step I),
and the following is a relatively straightforward consequence of 
a theorem of 
Breuil and Diamond \cite[Th\'{e}or\`{e}me 3.2.2]{BreuilDiamond}, which
in turn is a consequence of deep results due
to Kisin \cite{Kisin}, Gee \cite{Gee}, and 
Barnet-Lamb, Gee and Geraghty \cite{BGG1}, \cite{BGG2}.
\begin{thm}
\label{thm:modlift} 
Let $E$ be an elliptic curve over a totally real number field $K$,
and let $p\ne 2$ be a rational prime. Write $\overline{\rho}
=\overline{\rho}_{E,p}$. Suppose
\begin{enumerate}
\item[(i)] $\overline{\rho}$ is modular,
\item[(ii)] $\overline{\rho} (G_{K(\zeta_p)})$
is absolutely irreducible.
\end{enumerate}
Then $E$ is modular.
\end{thm}
The derivation of Theorem~\ref{thm:modlift} from the aforementioned
theorem of Breuil and Diamond is given in Section~\ref{sec:modlift}.
Sections~\ref{sec:modcurves} and~\ref{sec:images}
are devoted to background results on modular curves and
images $\overline{\rho}(G_{K(\zeta_p)})$ that are needed in the
remainder of the paper.
In Sections~\ref{sec:lt} and~\ref{sec:35} we explain how
Langlands--Tunnell (Step II) and
Wiles' 3--5 modularity switching argument (Step III) carry over
to the setting of elliptic curves over totally real fields.
These, together with Theorem~\ref{thm:modlift} give the following.
\begin{thm}\label{thm:35}
Let $p=3$ or $5$. Let $E$ be an elliptic curve over a totally real field $K$.
Suppose that $\overline{\rho}_{E,p}(G_{K(\zeta_p)})$ is absolutely irreducible.
Then $E$ is modular.
\end{thm}


We are interested
in proving modularity of elliptic curves over real quadratic fields. 
If we are to follow the same strategy as over $\Q$, we will 
have to examine the  non-cuspidal real quadratic points on certain modular
curves; these parametrize
elliptic curves where $\overline{\rho}_{E,3}(G_{K(\zeta_3)})$ and 
 $\overline{\rho}_{E,5}(G_{K(\zeta_5)})$ are simultaneously absolutely reducible. 
Some of these modular curves (for
example the elliptic curve $X_0(15)$) have infinitely many real quadratic points,
and there seems to be no simple way to prove that these points are modular.
We are forced to interpolate a further step
into the above strategy, by adding a 3--7 modularity switching argument (Section~\ref{sec:37}),
which incorporates some ideas of 
Manoharmayum \cite{Mano1}, \cite{Mano2}. 
\begin{thm}\label{thm:7}
Let $E$ be an elliptic curve over a totally real field $K$.
Suppose that $\overline{\rho}_{E,7}(G_{K(\zeta_7)})$ is absolutely irreducible.
Then $E$ is modular.
\end{thm}

Thanks to the addition of Theorem~\ref{thm:7},
 an elliptic curve $E$ over a totally real field $K$ 
is modular except possibly if the images $\overline{\rho}_{E,p}(G_{K(\zeta_p)})$
are simultaneously absolutely reducible for $p=3$, $5$, $7$.
In Section~\ref{sec:Karbproof} we show that such an elliptic curve gives
 rise to a $K$-point on one of $27$ modular curves:
\begin{equation}\label{eqn:27mod}
X(u,v,w) \qquad 
u \in \{\mathrm{b}3,\mathrm{s}3, \mathrm{ns}3\}, \quad
v \in \{\mathrm{b}5,\mathrm{s}5,\mathrm{ns}5\}, \quad
w \in \{\mathrm{b}7,\mathrm{s}7,\mathrm{ns}7\}.
\end{equation}
See Section~\ref{sec:modcurves} for an explanation of the
notation, but for now $\mathrm{b}$, $\mathrm{s}$ and $\mathrm{ns}$
respectively stand for \lq Borel\rq, \lq normalizer of split Cartan\rq\
and \lq normalizer of non-split Cartan\rq.
%
The curves in \eqref{eqn:27mod} all have genera $>1$, and so by Faltings' Theorem \cite{Faltings1}
possess only finitely many $K$-points. 
\begin{thm}\label{thm:Karb}
Let $K$ be a totally real field. All but finitely many $\overline{K}$-isomorphism
classes of elliptic curves over $K$ are modular. 
\end{thm}

Now a possible route to proving Theorem~\ref{thm:modquadratic}
is to attempt to enumerate all non-cuspidal real quadratic points
on the $27$ modular curves in~\eqref{eqn:27mod}
(or suitable quotients)
and then prove that the resulting points are modular.
This appears to be a monumental task,
and some of the curves on this list have
defeated us. 
In Sections~\ref{sec:imagesagain} and~\ref{sec:modcurvesagain}
we carry out a more precise examination
of cases where $\overline{\rho}_{E,p}(G_{K(\zeta_p)})$ 
is absolutely reducible, and the corresponding modular curves.
This allows us
 (Section~\ref{sec:modcurvesquad})  to
reduce the
proof of Theorem~\ref{thm:modquadratic} to 
establishing lemmas~\ref{lem:Qsqrt5} and~\ref{lem:7mod} below. 
\begin{lem}\label{lem:Qsqrt5}
Elliptic curves over $\Q(\sqrt{5})$ are modular.
\end{lem}
\begin{lem}\label{lem:7mod}
All non-cuspidal real quadratic points on the modular curves 
\begin{gather} \label{eqn:first3}
X(\mathrm{b}5,\mathrm{b}7), \qquad
X(\mathrm{b}3,\mathrm{s}5), \qquad
X(\mathrm{s}3,\mathrm{s}5),\\
\label{eqn:last4}
X(\mathrm{b}3,\mathrm{b}5,\mathrm{d}7), \quad
X(\mathrm{s}3,\mathrm{b}5,\mathrm{d}7), \quad
X(\mathrm{b}3,\mathrm{b}5,\mathrm{e}7), \quad
X(\mathrm{s}3,\mathrm{b}5,\mathrm{e}7),
\end{gather}
correspond to modular elliptic curves.
\end{lem}
The notation $\mathrm{d}7$ and $\mathrm{e}7$ is explained in
Section~\ref{sec:modcurvesagain}, but for now we remark
that they indicate mod $7$ level structures
that are respectively finer than \lq normalizer of split Cartan\rq\
and \lq normalizer of non-split Cartan\rq. 

The seven modular curves in \eqref{eqn:first3},~\eqref{eqn:last4} have genera $3$, $3$, $4$,
$97$, $153$, $73$ and $113$. It seems that these modular curves
are not amenable to 
the methods of Mazur \cite{MazurEisen}, Kamienny \cite{Kamienny} and 
Merel \cite{Merel} 
(we comment on this
at the end of Section~\ref{sec:modcurvesquad}). 
In the literature there are explicit methods which, in certain favourable circumstances,
are capable of determining all quadratic points on
a given curve (e.g.\ \cite{chabsym}, \cite{BN}), using an explicit model for the curve
and knowledge of the Mordell--Weil group of its Jacobian. As far as we are aware, such
methods have only been applied to curves of genus $2$ and $3$, and our computations
for the three curves in \eqref{eqn:first3} pushes these methods to their limit.
In Sections~\ref{sec:b5b7},~\ref{sec:b3s5},~\ref{sec:s3s5} we determine 
the quadratic points on the three modular curves in \eqref{eqn:first3},
and deduce modularity of the non-cuspidal real quadratic points.


The four modular curves in \eqref{eqn:last4} are extremely complicated on account of
their large genera, and certainly beyond the aforementioned explicit methods
found in the literature.
Instead we give a representation of these curves as normalizations
of fibre products of curves (Section~\ref{sec:models}), and develop an analogue of the
Mordell--Weil sieve that applies to such fibre products (Section~\ref{sec:sieve}). 
Remarkably, we are able to show (Section~\ref{sec:b3b5d7}) 
with the help of this sieve,
that the non-cuspidal real quadratic points on these curves have
$j$-invariants belonging to $\Q$, and so their modularity
follows from the modularity of elliptic curves over $\Q$.

The proof Theorem~\ref{thm:modquadratic}
is completed in Section~\ref{sec:Qsqrt5}
by proving Lemma~\ref{lem:Qsqrt5}.
This requires 
the computation of $\Q(\sqrt{5})$-points on another four modular curves.

We end the paper by 
briefly mentioning some further implications of this work 
for modularity of elliptic curves over arbitrary totally real fields
(Section~\ref{sec:totreal}),
and also possible Diophantine applications, including the Fermat equation
(Section~\ref{sec:dio}).

\medskip
 
\noindent \textbf{Acknowledgements.} 
We would like to thank the referees for useful comments.
It is a pleasure to express our sincere gratitude 
to a large number
of colleagues for their help and advice during the course
of writing this paper: 
Samuele Anni,
Alex Bartel,
Peter Bruin,
Frank Calegari,
Tommaso Centeleghe,
John Cremona, 
Lassina Demb\'{e}l\'{e},
Fred Diamond,
Luis Dieulefait,
Toby Gee,
Ariel Pacetti,
Richard Taylor
and 
Damiano Testa.
We would also like to thank
Rajender Adibhatla,
Shuvra Gupta,
Derek Holt, 
David Loeffler and
Panagiotis Tsaknias 
for useful discussions.

\medskip

\noindent \textbf{Computational Details.} 
Some of the computational
details in Sections~\ref{sec:b5b7}--\ref{sec:Qsqrt5}
were carried out in 
the computer algebra system {\tt Magma} \cite{MAGMA}.
The reader can find the {\tt Magma} scripts 
for verifying these \label{magma}
computations 
(together with extensive commentary)
at:\newline
{\tt http://arxiv.org/abs/1310.7088} 

Our programs make use of several {\tt Magma} packages. We would like to acknowledge
these packages and their authors and main contributers:
\begin{itemize}
\item {\tt Hyperelliptic Curves}: Nils Bruin, Brendan Creutz, Steve Donnelly, Michael Harrison, David Kohel, Michael Stoll, 
Paul van Wamelen;
\item {\tt Small Modular Curves Database}: Michael Harrison;
\item {\tt Algebraic Function Fields}: Florian He{\ss}, Claus Fieker, Nicole Sutherland;
\item {\tt Modular Forms}: William Stein, Kevin Buzzard, Steve Donnelly;
\item {\tt Modular Abelian Varieties}: William Stein, Jordi Quer;
\item {\tt Matrix Groups over Finite Fields}: Eamonn O'Brien.
\end{itemize}

\part{Background}

\section{Modularity Lifting: Proof of Theorem~\ref{thm:modlift}}\label{sec:modlift}
We assume the hypotheses of Theorem~\ref{thm:modlift}.
Thus $K$ is a totally real field, $E$ is 
an elliptic curve over $K$, and $p \ne 2$ is
a rational prime so that $\overline{\rho}:=\overline{\rho}_{E,p}$
is modular, and $\overline{\rho}(G_{K(\zeta_p)})$
is  absolutely irreducible.

Write $\rho : G_K \rightarrow \GL_2(\Z_p)$
for the representation arising from the action of $G_K$ 
on the Tate module $T_p(E)$.
For any finite place $\upsilon$ of $K$, we write $D_\upsilon \subset G_K$
for the corresponding decomposition group. 
We apply Th\'{e}or\`eme 3.2.2 of \cite{BreuilDiamond}
to $\overline{\rho}$. That theorem involves making certain choices
and we merely make the obvious ones as follows.
We let $S$ be the set of finite places $\upsilon$ of $K$ of bad reduction for $E$
together with the places above $p$. We let $T$ be the
subset of places above $p$ where $E$ has potentially
ordinary reduction or potentially multiplicative reduction. 
In the notation of that theorem $\epsilon : G_K \rightarrow \Z_p^*$
is the $p$-adic cyclotomic character, and we choose $\psi=\epsilon$.
 We let $\rho_\upsilon=\rho \vert_{D_\upsilon}$.
For a place $\upsilon$ we write $[r_\upsilon,N_\upsilon]$ for the Weil--Deligne
type of $\rho_\upsilon$. We note
(\cite[Sections 14 and 15]{Rohrlich2}) that $N_\upsilon=0$ if and only if
$E$ has potentially good reduction at $\upsilon$. If $p\ne 5$, the conditions of \cite[Th\'{e}or\`{e}me 3.2.2]{BreuilDiamond}
are satisfied; that theorem asserts existence of
a global modular lift of $\overline{\rho}$ having the same local
properties as $\rho$ (properties (i)--(v)
in the conclusion of the theorem). It also asserts that every global lift
having these properties 
is modular, so $\rho$ is modular. When $p=5$, the theorem has an extra hypothesis: namely that
the projective image of $\overline{\rho}\vert_{G_{K(\zeta_5)}}$ is not
$\PSL_2(\F_5)$. The reason for this additional hypothesis is to assure the existence of Taylor--Wiles systems
(as in \cite[Section 3.2.3]{Kisin}). However, in our situation, because $\det\overline{\rho}$
is the mod $5$ cyclotomic character, we can still choose Taylor--Wiles systems
without this extra condition, as we now explain. In the notation of \cite[Section 3.2.3]{Kisin} (except
we change the name of the field to $K$) and following the proof of Theorem 2.49
in \cite{DDT}, what we need to show is that for each positive integer $n$ and 
non-trivial cocycle $\psi\in \mathrm{H}^1(G_{K,S},\textrm{ad}^0 \overline{\rho}(1))$, we
can find a place $\upsilon\notin S$ of $K$ such that 
\begin{enumerate}
\item[(i)] $\lvert k(\upsilon) \rvert \equiv 1 \pmod{p^n}$, and $\overline{\rho}(\mathrm{Frob}_\upsilon)$ has distinct eigenvalues.  
\item[(ii)]   
The image of $\psi$ under the restriction map 
\begin{align*}
  \mathrm{H}^1(G_K,\textrm{ad}^0 \overline{\rho}(1))\to \mathrm{H}^1(D_\upsilon,\textrm{ad}^0 \overline{\rho}(1))
\end{align*}
is non-trivial.
\end{enumerate}
 Let $K_m$ be the extension of $K(\zeta_{p^m})$ cut out by
$\textrm{ad}^0\overline{\rho}$. The argument in \cite[Theorem 2.49]{DDT} then
works once we know the restriction of $\psi$ to
$\mathrm{H}^1(G_{K_0},\textrm{ad}^0\overline{\rho}(1))$ is non-trivial. To do this we
want to show that
$\mathrm{H}^1(\Gal(K_0/K),\textrm{ad}^0\overline{\rho}(1)^{G_{K_0}})=0$. Because
$G_{K_0}$ acts trivially on $\textrm{ad}^0\overline{\rho}$, the coefficient
module vanishes unless $\zeta_p\in K_0$. Because
$\mathrm{H}^1(\PSL_2(\F_5),\Sym^2\F_5^2)=\F_5$ does not vanish, the extra
condition when $p=5$ was used to exclude the possibility that
$\Gal(K_0/K)=\PSL_2(\F_5)$ and $\zeta_5\in K_0$. But using the fact that
$\overline{\rho}$ has cyclotomic determinant, we can get rid of it as follows.
Suppose the above happens, then $\sqrt{5}\in K$ because
$\det:\overline{\rho}(G_K) \to \PGL_2(\F_5)\to \F_5^*/(\F_5^*)^2$ is induced by
the cyclotomic character. But $\zeta_5\in K_0\setminus K$ implies that
$\Gal(K_0/K)$ 
has a quotient of order $2$. 
As
$\Gal(K_0/K)\cong \PSL_2(\F_5) \cong A_5$, we have a contradiction.

\section{Background on Modular Curves}\label{sec:modcurves}
In order to apply Theorem~\ref{thm:modlift}, we need a better understanding
of its condition (ii). 
We will show that if condition (ii) fails, then the 
image of $\overline{\rho}_{E,p}$ lies 
either in a Borel subgroup of $\GL_2(\F_p)$, or in certain subgroups
contained in the normalizers of Cartan subgroups. In this case,
the elliptic 
curve $E/K$ gives rise to a non-cuspidal $K$-point on a
modular curve corresponding to one of these subgroups. We now
briefly sketch the relationship between certain subgroups
of $\GL_2(\F_p)$ and the corresponding modular curves.
For more details we recommend the articles by 
Mazur \cite{Mazur} and by Rohrlich \cite{Rohrlich}.

Let $p$ be a prime, and $H$ a subgroup of $\GL_2(\F_p)$ 
satisfying the condition $-I \in H$. Write
$H_0=H \cap \SL_2(\F_p)$ and denote by $\widetilde{H}$
the preimage of $H_0$ via the natural map
$\SL_2(\Z) \rightarrow \SL_2(\F_p)$.
Then $\tilde{H}$ is a congruence subgroup of $\SL_2(\Z)$ of level $p$
and we define the modular curve
corresponding to $H$ by 
\[
X(H):=\widetilde{H} \backslash
\mathfrak{H}^*, \qquad \mathfrak{H}^*:=\mathfrak{H} \cup \Q \cup \{\infty\}
\]
where $\mathfrak{H}$ is upper half-plane.  
The above construction gives $X(H)$ as compact Riemann surface, or 
equivalently as an algebraic curve over $\C$. It turns out however
that $X(H)$ is defined over the field $\Q(\zeta_p)^{\det(H)}$;
here we view $\det(H)$ as subgroup of $\Gal(\Q(\zeta_p)/\Q)$
via the cyclotomic character $\chi_p : \Gal(\Q(\zeta_p)/\Q)
\xrightarrow{\cong} \F_p^*$. We shall assume henceforth that 
$\det(H)=\F_p^*$, and so $X(H)$ is defined over $\Q$. The curve
$X(H)$ admits a natural map 
\[
j: X(H) \rightarrow X(1), \qquad  (X(1):=\SL_2(\Z) \backslash \mathfrak{H}^*
\cong \PP^1)
\]
corresponding to the inclusion $\tilde{H} \subset \SL_2(\Z)$.
The map $j$ is defined over $\Q$, and is ramified only above
$0$, $1728$ and $\infty$. We refer to the points of the fibre
$j^{-1} (\infty) \subset X(H)(\overline{\Q})$ as the cusps
of $X(H)$. For a number field $K$, the non-cuspidal
$K$-points of $X(H)$ correspond to $\overline{\Q}$-isomorphism classes of pairs $(E,\eta\circ H)$ of elliptic curves
$E/K$ with an $H$-orbit of isomorphisms 
\begin{align*}
  \eta: \F_p^2
\xlongrightarrow{\cong}
E[p](\overline{\Q})
\end{align*}
such that $\eta^{-1} \sigma \eta\in H$ for all $\sigma\in G_K$ (note that $H$ acts on the set of such $\eta$ by its action on the domain of $\eta$). An elliptic
curve $E$ over $K$ supports a non-cuspidal point on $X(H)$ if and only if the
image $\overline{\rho}_{E,p}(G_K)$ is conjugate to a subgroup of $H$. In
general, it is possible for one $E$ to support distinct points of $X(H)$. In
fact, picking an $\eta$, the number of points
supported on $E$ is the number of $G_K$-fixed points of the double coset space
\[
\eta^{-1}\cdot \Aut_{\overline{\Q}}(E) \cdot \eta \, \backslash \, \GL_2(\F_p) \, /\, H \, ,
\]
 where the action
of $G_K$ depends on the chosen $\eta$. Moreover if $E/K$ gives rise to the $K$-point $P$
on  $X(H)$ then $j(P)$ is the $j$-invariant of $E$. If two elliptic curves
$E_1$, $E_2$ defined
over $K$, and having $j$-invariants $\ne 0$, $1728$ give rise to the same
$K$-point of $X(H)$, then they have the same $j$-invariant and hence are
quadratic twists of each other.

\subsection{Borel Subgroups, Cartan Subgroups, and Normalizers}

Let $p\ge 3$ be a prime.
Let $B(p)$ denote the Borel subgroup of $\GL_2(\F_p)$.
Let $C_\mathrm{s}(p)$ and $C_\mathrm{ns}(p)$ respectively
denote the split and non-split Cartan subgroups of $\GL_2(\F_p)$,
and let $C_\mathrm{s}^+(p)$ and $C_\mathrm{ns}^+(p)$ 
respectively be their normalizers. The three groups
$B(p)$, $C_\mathrm{s}^+(p)$ and $C_\mathrm{ns}^+(p)$
contain $-I$ and have determinants equal to $\F_p^*$.
Thus they correspond to modular curves defined over $\Q$.
The usual notation for these modular curves is
$X_0(p)$, $X_\mathrm{split}(p)$ and $X_\mathrm{nonsplit}(p)$.
We shall find it convenient to denote these modular curves
by the following non-standard notation: $X(\mathrm{b}p)$, $X(\mathrm{s}p)$ and $X(\mathrm{ns}p)$,
so for example $X_0(5)$, $X_\mathrm{split}(5)$ and $X_\mathrm{nonsplit}(5)$
become $X(\mathrm{b}5)$, $X(\mathrm{s}5)$, $X(\mathrm{ns}5)$.

\subsection{Fibre Products}
Now let $p_1,\dots,p_r \ge 3$ be distinct primes. 
Let $H_i \subseteq \GL_2(\F_{p_i})$ be subgroups,
all satisfying $-I \in H_i$ and $\det(H_i)=\F_{p_i}^*$.
Let $\widetilde{H_i}$ be the corresponding congruence subgroups of $\SL_2(\Z)$ defined above.
Then $\bigcap \widetilde{H_i}$ is a congruence subgroup of $\SL_2(\Z)$ of level $\prod p_i$.
We define $X(H_1,\ldots,H_r)$ to be the modular curve
\[
X(H_1,\ldots,H_r):=\left(\widetilde{H_1}\cap \dots \cap \widetilde{H_r}\right) \, \backslash \,
\mathfrak{H}^*.
\]
This again is defined over $\Q$, and is the normalization of 
the fibre product
\[
X(H_1) \times_{X(1)} X(H_2) \times_{X(1)} \cdots \times_{X(1)} X(H_r).
\]
Again 
$X(H_1,\dots,H_r)$ is naturally endowed with a $j$-map to $X(1)$
coming from the inclusion $\widetilde{H_1}\cap \dots \cap \widetilde{H_r}
\subseteq \SL_2(\Z)$. This agrees
with the composition of the projection onto
any of the $X(H_i)$ with $j: X(H_i) \rightarrow X(1)$.
If $K$ is a number field, 
non-cuspidal $K$-points correspond to (isomorphism classes of) pairs $(E,\eta\circ \prod H_i)$ of elliptic curves
$E/K$ with an $\prod H_i$-orbit of isomorphisms 
\begin{align*}
  \eta:\prod \F_{p_i}^2 
\xlongrightarrow{\cong}
E[\prod p_i](\overline{\Q})
\end{align*}
such that $\eta^{-1}\overline{\rho}_{E,p_i}(G_K)\eta$ is
contained in $H_i$, for $i=1,\ldots,r$. 

In the notation $X(H_1,\ldots,H_r)$, we shall replace
the group $H_i$ with 
$\mathrm{b}p$,
$\mathrm{s}p$, $\mathrm{ns}p$
if it equals 
$B_0(p)$, $C_\mathrm{s}^+(p)$
or $C_\mathrm{ns}^+(p)$
respectively.
Thus for example,
$X(\mathrm{b}5,\mathrm{b}7)$
is isomorphic to $X(\mathrm{b}5) \times_{X(1)} X(\mathrm{b}7)$ 
and in standard notation can be denoted by $X_0(35)$.

\section{The Image $\overline{\rho}_{E,p}(G_{K(\zeta_p)})$}\label{sec:images}
The following proposition is our first attempt at understanding
condition (ii) of Theorem~\ref{thm:modlift}.
It relates the failure
of condition (ii) to $K$-points on the modular curves $X(\mathrm{b}p)$, $X(\mathrm{s}p)$
and $X(\mathrm{ns}p)$. 

\begin{prop}\label{prop:large}
Let $E$ be an elliptic curve over a totally real field $K$. Let $p \ge 3$ be
a rational prime, and write $\overline{\rho}=\overline{\rho}_{E,p}$.
Then
\begin{enumerate}
\item[(i)] $\overline{\rho}(G_{K(\zeta_p)})=\overline{\rho}(G_K) \cap \SL_2(\F_p)$.
\item[(ii)] If $\overline{\rho}(G_{K(\zeta_p)})$ is absolutely reducible,
then $\overline{\rho}(G_K)$ is contained either in a Borel subgroup,
or in the normalizer of a Cartan subgroup. In this case $E$ gives
rise to a non-cuspidal $K$-point on $X(\mathrm{b}p)$, $X(\mathrm{s}p)$
or $X(\mathrm{ns}p)$.
\end{enumerate}
\end{prop}
\begin{proof}
The determinant of $\overline{\rho}$ is the mod $p$ cyclotomic character 
\[
\chi_p : G_K \rightarrow \Aut(\mu_p) \cong \F_p^*,
\]
where $\mu_p$ is the group of $p$-th roots of unity.
The first part follows immediately.

Write $G=\overline{\rho}(G_K)$, and suppose $G$
is not contained in a Borel subgroup. Then $G$ is irreducible,
and as $K$ is totally real it is absolutely irreducible.
The second part follows from Lemma~\ref{lem:exceptional}
below. 
\end{proof}

\begin{lem}\label{lem:exceptional}
Let $p \ge 3$ be prime. Let $G$ be an absolutely irreducible subgroup 
of $\GL_2(\F_p)$ such that $G \cap \SL_2(\F_p)$ is absolutely
reducible. Then
$G$ is contained in the normalizer of a Cartan subgroup.

If moreover $\det(G)=\F_p^*$,
then $G \cap \GL_2^+(\F_p)$ is contained in a Cartan subgroup, where
$\GL_2^+(\F_p)$ is the subgroup of $\GL_2(\F_p)$
consisting of matrices with square determinant.
\end{lem}
\begin{proof}
Write $G^\prime=G \cap \SL_2(\F_p)$.
If $p$ divides the order of $G$ then 
(e.g.\ \cite[Section 2]{Serre}
or \cite[Lemma 2]{SwD}) 
 either $G$ contains $\SL_2(\F_p)$ or $G$ is 
contained in a Borel subgroup. In the former case, $G^\prime$
contains $\SL_2(\F_p)$ and so is absolutely irreducible,
and in the latter case, $G$ is reducible. We may thus
suppose that $p$ does not divide the order of $G$.

Write $H$ for the image of $G$ in $\PGL_2(\F_p)$, and let $n$
be its order.
Let $\cM$
be the set of points $P \in \PP^1(\overline{\F}_p)$
that are fixed by at least one non-trivial element of $H$.
Swinnerton-Dyer shows \cite[pages 13--14]{SwD} that $\cM$ is 
stable under the action of $G$,
and that it can be written as the union of 
disjoint orbits $\cN_1,\dots,\cN_\nu$, with the following 
exhaustive list of possibilities.
\begin{enumerate}
\item[(a)] 
$H$ is cyclic, $G$ is contained 
in a Cartan subgroup:  
$\nu=2$,  $\lvert \cN_1 \rvert=
\lvert \cN_2 \rvert=1$;
\item[(b)] 
$H$ is dihedral, $G$ is contained in 
the normalizer of a Cartan subgroup: 
$\nu=3$, $n\ge 4$ even, $\lvert \cN_1 \rvert=\lvert \cN_2
\rvert=n/2$, $\lvert\cN_3 \rvert =2$;
\item[(c)] 
$H$ is isomorphic to $A_4$: 
$\nu=3$, $n=12$, $\lvert \cN_1 \rvert=6$, $\lvert \cN_2 \rvert=
\lvert \cN_3 \rvert=4$;
\item[(d)] 
$H$ is isomorphic to $S_4$:
$\nu=3$, $n=24$, $\lvert \cN_1 \rvert=12$, $\lvert \cN_2 \rvert=8$, 
$\lvert \cN_3 \rvert =6$;  
\item[(e)] 
$H$ is isomorphic to $A_5$:
$\nu=3$, $n=60$, $\lvert \cN_1 \rvert=30$, $\lvert \cN_2 \rvert=20$, 
$\lvert \cN_3 \rvert=12$.  
\end{enumerate}

In our situation, case (a) is easy to rule out as it forces
$G$ to be absolutely reducible---indeed, the unique
elements of $\cN_1$, $\cN_2$ are fixed by $G$.
We would like to rule out cases (c), (d), (e).

Observe that $G^\prime$ is the kernel of the restriction 
to $G$ of $\det : G \rightarrow \F_p^*$. Thus $G^\prime$
is normal in $G$ with cyclic quotient. If $G^\prime$
consists only of scalar matrices, then $H$ is
a quotient of $G/G^\prime$ and so is cyclic, which
contradicts possibilities (b)--(e). Thus $G^\prime$ contains
some non-scalar matrix $g$. The image of $g$ in $H$
is non-trivial.

As $G^\prime$ is absolutely reducible,
there is some  $P \in \PP^1(\overline{\F}_p)$ which is
 fixed by $G^\prime$. In particular, $P$ belongs to 
one of the orbits $\cN_i$, which we denote by $\cN$.
 If $Q$ is some other element of  $\cN$ then its stabilizer
in $G$ is conjugate to the stabilizer of $P$.
However, $G^\prime$ is normal in $G$ and is contained in the stabilizer
of $P$. Thus $G^\prime$ stabilizes $Q$ too. It follows that every
element $Q$ of $\cN$ is an eigenvector of $g$.
However, $g$ is a non-scalar and $p$ does not divide
its order (as $g \in G$). Thus $g$ has precisely two
eigenvectors in $\PP^1(\overline{\F}_p)$.
Hence $\lvert \cN \rvert \le 2$.
From the above list of cases (b)--(e), we see that we are in case (b),
and $\lvert \cN \rvert=2$. This proves the first part of the lemma.

For the second part, suppose that
$\det(G)=\F_p^*$. 
Let $G^+$ be the stabilizer of $P$, which by the above contains $G^\prime$.
Since $G$ acts transitively on $\cN$,
which has size $2$, the subgroup $G^+$ must have index $2$. 
From the  isomorphism $\det : G /G^\prime \rightarrow \F_p^*$
we see that $\det(G^+/G)$ consists of the squares in $\F_p^*$.
Hence $G^+=G \cap \GL_2^+(\F_p)$. Now $G^+$ stabilizes $P$,
and by the argument above also stabilizes the other element of $\cN$. Thus 
$G^+$ is contained in a Cartan subgroup.
\end{proof}

\section{Langlands--Tunnell: Proof of Theorem~\ref{thm:35} for $p=3$}\label{sec:lt}
We shall need the following consequence of the well-known Langlands--Tunnell
Theorem.
\begin{thm}[Langlands--Tunnell]
\label{thm:LT} 
Let $E$ be an elliptic curve over a totally real number field $K$
and suppose $\overline{\rho}_{E,3}$ is irreducible. Then 
$\overline{\rho}_{E,3}$ is modular.
\end{thm}
The deduction of Theorem~\ref{thm:LT} from the Langlands--Tunnell
Theorem \cite{Langlands}, \cite{Tunnell} is well-known; a proof is given 
by 
Dieulefait and Freitas in \cite[Lemma 4.2]{DF2}
and follows the lines of the argument by Wiles \cite[Chapter 5]{Wiles}
for classical modular forms.

The following is an immediate corollary of Langlands--Tunnell
and Theorem~\ref{thm:modlift}. It is of course Theorem~\ref{thm:35} with $p=3$.
\begin{cor}\label{cor:3mod}
If $E$ is an elliptic curve over a totally real field $K$
and $\overline{\rho}_{E,3}(G_{K(\zeta_3)})$ is absolutely irreducible,
then $E$ is modular.
\end{cor}

\part{Modularity Switching}

\section{3--5 Modularity Switching: Proof of Theorem~\ref{thm:35} for $p=5$}
\label{sec:35}
In this section we prove Theorem~\ref{thm:35} for $p=5$.
This proof follows from Wiles' original 3--5 modularity switching argument
over $\Q$ \cite[Proof of Theorem 5.2]{Wiles}. For the convenience of the reader
and as a warm-up for the more elaborate switching argument in the next section,
we include a sketch. 
\begin{lem}\label{lem:xe5}
Let $E$ be an elliptic curve over a number field $K$.
Then there is an elliptic curve $E^\prime/K$ such that
\begin{enumerate}
\item[(i)] $\overline{\rho}_{E,5} \sim \overline{\rho}_{E^\prime,5}$;
\item[(ii)] $\overline{\rho}_{E^\prime,3}(G_K)$ contains $\SL_2(\F_3)$.
\end{enumerate}
\end{lem}
We first show how the lemma implies Theorem~\ref{thm:35} for $p=5$. 
\begin{proof}[Proof of Theorem~\ref{thm:35} for $p=5$]
Let $E$ be an elliptic curve over a totally real number field $K$ and let $E^\prime$
be as in the statement of Lemma~\ref{lem:xe5}. 
It follows from 
Corollary~\ref{cor:3mod}
that $E^\prime$ is modular. Thus $\overline{\rho}_{E,5}\sim \overline{\rho}_{E^\prime,5}$
is modular. Applying Theorem~\ref{thm:modlift} to $E$
with $p=5$, 
we see that $E$ is modular. 
\end{proof}
Let us now prove lemma \ref{lem:xe5}. Let $K$ be a number field and $E/K$ is an
elliptic curve. Let us denote by $X_E(5)$ the twisted modular curve over $K$
which classifies elliptic curves $E^\prime$ with a symplectic isomorphism of
group schemes $E[5] \rightarrow E^\prime[5]$. Here the adjective {\bf
symplectic} means that the isomorphism is compatible with the Weil pairings on
$E[5]$ and $E^\prime[5]$. We observe the following properties
\begin{itemize}
\item  $X_E(5)\cong \mathbb{P}^1$ over $K$. This is because $X_E(5)$ isomorphic
(over $\overline{K}$)
to the standard modular curve $X(5)$, 
and hence is a smooth curve of genus $0$ over
$K$, possessing by construction a $K$-rational point (coming from $E$
itself).
\item 
There is bijection between the non-cuspidal $K$-points on $X_E(5)$ and classes
of pairs
$(E^\prime,u)$ where $E^\prime$ is an elliptic curve over 
$K$ and $u$ is a symplectic isomorphism $E[5] \rightarrow E^\prime[5]$
of $G_K$-modules. Two such pairs $(E_1^\prime,u_1)$ and $(E_2^{\prime}, u_2)$
belong to the same class if there is a $K$-isomorphism $\phi:E_1^\prime \rightarrow E_2^\prime$
such that $u_2=\phi \circ u_1$.
\end{itemize}
Thus a non-cuspidal $K$-point of $X_E(5)$ gives an elliptic curve $E^\prime$
satisfying condition (i), and we want to find such a point where (ii) holds.
The maximal subgroups of $\GL_2(\F_3)$ are $\SL_2(\F_3)$, the Borel subgroups,
and the normalizers of non-split Cartan subgroups~\footnote{In $\GL_2(\F_3)$
the normalizer of a split Cartan subgroup is contained
as an index $2$ subgroup in the normalizer of a non-split Cartan subgroup, as
the latter is a 2-Sylow subgroup.}.
If the image of  $\overline{\rho}_{E^\prime, 3}$ does not contain
$\SL_2(\F_3)$ then $E^\prime$ gives rise to a $K$-point of $X_E(5)\times_{X(1)}
X(\mathrm{b}3)$ or $X_E(5)\times_{X(1)} X(\mathrm{ns}3)$. Over $\mathbb{C}$,
the normalizations of these two curves are isomorphic to the modular curves
$X(5,\mathrm{b}3)$ and $X(5,\mathrm{ns}3)$, and hence are irreducible. It
follows from Hilbert's irreducibility theorem \cite[Theorem 3.4.1]{SeGalois}
that there is a non-cuspidal $K$-rational point of $X_E(5)$ which does not
lift to a $K$-rational point of  $X_E(5)\times_{X(1)} X(\mathrm{b}3)$,
$X_E(5)\times_{X(1)} X(\mathrm{ns}3)$. This point provides the required
$E^\prime$.

\section{3--7 Modularity Switching: Proof of Theorem~\ref{thm:7}}
\label{sec:37}
In this section we prove Theorem~\ref{thm:7}.
The proof resembles Wiles' 3--5 
modularity
switching argument 
sketched in the previous section.
In fact, similar arguments 
 have been used to prove special cases
of Serre's modularity conjecture for $\GL_2(\F_q)$ for $q=5$, $7$ and $9$
by Taylor \cite{TaylorIco}, 
Manoharmayum \cite{Mano1}, \cite{Mano2} 
and Ellenberg \cite{ellenberg}.
Our argument in this section is closest in spirit to Manoharmayum's;
however we benefit from having at our disposal a more powerful 
modularity lifting theorem, allowing us to work in greater generality.
We shall need the following lemma.

\begin{lem}\label{lem:xe7}
Let $E$ be an elliptic curve over a totally real number field $K$.
Then there is an extension $L/K$ and an elliptic curve
$E^\prime/L$ such that
\begin{enumerate}
\item[(i)] $L$ is totally real;
\item[(ii)] $[L:K]=4$;
\item[(iii)] $L \cap K(E[7])=K$;
\item[(iv)] $\overline{\rho}_{E,7} \vert_{G_L} \sim \overline{\rho}_{E^\prime,7}$;
\item[(v)] $\overline{\rho}_{E^\prime,3}(G_L)$ contains $\SL_2(\F_3)$.
\end{enumerate}
\end{lem}

Before proving Lemma~\ref{lem:xe7} we use it prove Theorem~\ref{thm:7}.
\begin{proof}[Proof of Theorem~\ref{thm:7}]
Let $E$ be an elliptic curve over a totally real number field $K$ and let $L$ and $E^\prime$
be as in the statement of Lemma~\ref{lem:xe7}. It follows 
from Corollary~\ref{cor:3mod} that
$E^\prime$ is modular. Thus $\overline{\rho}_{E^\prime,7}$
is modular. Hence $\overline{\rho}_{E,7} \vert_{G_L}$ is modular. 

From now on, we shall only be concerned with the elliptic
curve $E$, and will write $\overline{\rho}=\overline{\rho}_{E,7}$.
We will shortly show that $\overline{\rho} (G_{L(\zeta_7)})$
is absolutely irreducible. 
Thus we may apply Theorem~\ref{thm:modlift} to deduce that $E/L$ is modular.
Now by a repeated application of cyclic base change results of Langlands \cite{Langlands},
it follows that $E/K$ is modular.

All that remains to show is that $\overline{\rho}_{E,7} (G_{L(\zeta_7)})$
is absolutely irreducible. It is here that we make use of the hypothesis in Theorem~\ref{thm:7}
that
$\overline{\rho} (G_{K(\zeta_7)})$ is absolutely irreducible. 
By Proposition~\ref{prop:large},
\[
\overline{\rho}(G_{K(\zeta_7)})=\overline{\rho}(G_K) \cap \SL_2(\F_7),
\qquad 
\overline{\rho}(G_{L(\zeta_7)})=\overline{\rho}(G_L) \cap \SL_2(\F_7).
\]
To complete the proof it is enough to show that $\overline{\rho}(G_K)=\overline{\rho}(G_L)$.
The former is the Galois group of $K(E[7])/K$ and the latter is the Galois group
of $L(E[7])/L$. Since $L \cap K(E[7])=K$ by part (iii) of Lemma~\ref{lem:xe7} we see
that these groups are the same.
\end{proof}

We now turn our attention to proving Lemma~\ref{lem:xe7}.
As before $K$ is a totally real number field and $E/K$ is
an elliptic curve. We shall need to work with the modular curve $X_E(7)$.
For detailed expositions of the construction and properties of $X_E(7)$
we recommend the paper of 
Halberstadt and Kraus \cite{HK},  and that of Poonen, Schaefer and Stoll \cite{PSS}.
For our needs, the following basic facts found in the above references are sufficient.
\begin{itemize}
\item $X_E(7)$ is a smooth plane quartic curve defined over $K$.
\item $X_E(7)$ is a twist of the Klein quartic
\begin{equation}\label{eqn:klein}
X(7)\; :\; x^3 y+ y^3 z + z^3 x=0.
\end{equation}
\item For any extension $L/K$, there is bijection between the
non-cuspidal $L$-points on $X_E(7)$ and classes of pairs
$(E^\prime,u)$ where $E^\prime$ is an elliptic curve over 
$L$ and $u$ is a symplectic isomorphism $E[7] \rightarrow E^\prime[7]$
of $G_L$-modules.
Two such pairs $(E_1^\prime,u_1)$ and $(E_2^{\prime}, u_2)$
belong to the same class if there is an $L$-isomorphism $\phi:E_1^\prime \rightarrow E_2^\prime$
such that $u_2=\phi \circ u_1$.
\end{itemize}

An elliptic curve $E^\prime$ as in the statement of Lemma~\ref{lem:xe7}
corresponds to a point on $X_E(7)$ over a certain
totally real quartic extension $L/K$ satisfying
some additional hypotheses.
We shall generate points
on $X_E(7)$ defined over quartic extensions by 
taking the intersection of $X_E(7)$ with 
suitable $K$-lines. 
Write $\pdual$ for the dual space of $\PP^2$---this
is the parameter space of lines in $\PP^2$, and is
in fact isomorphic to $\PP^2$. To generate a point on
$X_E(7)$ that satisfies (ii)--(v) we
will choose a line $\ell \in \pdual(K)$ that
avoids certain thin subsets of $\pdual(K)$. To satisfy (i),
in addition to avoiding these thin subsets, we want our line
to belong to a certain open subset of $\pdual(K_\sigma)$
for each embedding $\sigma : K \hookrightarrow \R$. 

It is helpful to recall the notion of thin subsets 
\cite[Chapter 9]{SeMordell}, \cite[Chapter 3]{SeGalois}. Let $n \ge 1$.
By a {\bf thin subset of  $\PP^n(K)$ of type I} we 
mean a subset contained in $V(K)$ where $V$ is 
a proper Zariski-closed $K$-subvariety of $\PP^n$.
By a {\bf thin subset of $\PP^n(K)$ of type II}
we mean a subset contained in 
$\pi (V(K))$ where
$\pi: V \rightarrow \PP^n$ is a generically
surjective $K$-morphism of degree $\ge 2$
from an irreducible variety $V$ of dimension $\dim(V)=n$.
A {\bf thin subset of $\PP^n(K)$} is a subset
contained in the finite union of 
thin subsets of types I and II. Hilbert's Irreducibility
Theorem \cite[Theorem 3.4.1]{SeGalois} asserts that  a thin 
subset cannot equal $\PP^n(K)$.

\begin{proof}[Proof of Lemma~\ref{lem:xe7}] 
Fix a totally real field $K$
and an elliptic curve $E/K$. Write $X=X_E(7) \subset \PP^2$.
We shall write $X^{(n)}$ for the $n$-th symmetric power
of $X$ (see for example \cite[Section 3]{Milne}). This is a smooth $n$-dimensional variety over $K$, whose
$K$-points correspond to effective degree $n$ divisors on $X$
that are stable under the action of $G_K$.
Let $\phi$ be the morphism
\[
\phi: \pdual \hookrightarrow X^{(4)}, \qquad \ell \mapsto \ell \cdot X.
\]
We define certain thin subsets $S_1,\dots, S_{6}$ of $\pdual(K)$
and choose the line to lie outside these subsets. 
We let $S_1$ be the set of lines $\ell \in \pdual(K)$ passing through
any point belonging to either of these two finite sets:
\begin{enumerate}
\item[(a)] $\overline{K}$-cusps of $X_E$;
\item[(b)] $X_E(K(E[7]))$ (which is certainly a finite set by Faltings' Theorem \cite{Faltings1}).
\end{enumerate}
The set $S_1$ is plainly a union of thin subsets of type I.
Let $S_2$ be the set of lines $\ell \in \pdual(K)$ that
are tangent to $X$. These belong to the dual curve of $X$,
and hence $S_2$ is a thin subset of type I. 
Now any line belonging to $\pdual(K) \setminus S_1 \cup S_2$ intersects
$X$ either in an irreducible divisor of degree $4$, or in the sum
of two irreducible divisors of degree $2$. 

Write $\epsilon$ for
the natural map 
\[
\epsilon \; : \; X^{(2)} \times X^{(2)} \rightarrow X^{(4)}, \qquad (D_1,D_2) \mapsto D_1+D_2. 
\]
We let $S_3$ be the set of lines $\ell \in \pdual(K)$ so that $\phi(\ell)$ belongs to
$\epsilon (X^{(2)}(K) \times X^{(2)}(K))$. The set  $S_3$ certainly contains all $K$-lines $\ell$
so that $\ell \cdot X$ is the sum of two irreducible divisors of degree $2$. We would like to show that
$S_3$ is a thin set.
It is known for a curve $C$ of genus at least $2$
that $C^{(2)}(K)$ is finite, unless $C$ is hyperelliptic or bielliptic \cite[Corollary 3]{HS}.
In our situation $X$ is not hyperelliptic, but might be bielliptic \cite[Section 4]{bars2}, \cite[Th\'eor\`eme 8.1]{HK}. 
At any rate, $X^{(2)}$ is isomorphic to its image in the Jacobian of $X$ under an Abel-Jacobi map.
It follows from a theorem of Faltings \cite{Faltings} that
\begin{equation}\label{eqn:x2}
X^{(2)}(K)=\{D_1,\dots, D_n \} \cup \{ E_1(K), \dots, E_m(K) \}
\end{equation}
where the $D_i$ are finitely many divisors of degree $2$ and the $E_i$ are finitely many
elliptic curves defined over $K$ lying on $X^{(2)}$. Now
$\epsilon(X^{(2)}(K) \times X^{(2)}(K))$ consists of a union
of isolated points $D_i+D_j$, the $K$-points of curves $D_i \times E_j$ and  
the $K$-points of surfaces $E_i \times E_j$. The former two
categories will contribute to $S_3$ thin sets of type I. 
To show that $S_3$ is a thin set, it is enough to show for any $i$, $j$ that there is no 
$K$-rational map $\psi : \pdual \dashedrightarrow E_i \times E_j$,
birational onto its image, making the following diagram commute:
\[
\xymatrixcolsep{5pc}
\xymatrix{
& E_i \times E_j \ar@{^{(}->}[d] \\
& X^{(2)} \times X^{(2)}
 \ar^{\epsilon}[d]\\
\pdual \ar@{-->}^{\psi}[ruu]
\ar@{^{(}->}_{\phi}[r]
& X^{(4)} \\
}
\]
Any $E_i \times E_j$ is an abelian
variety and the image of any rational map $\pdual \dashedrightarrow E_i \times E_j$
must be constant \cite[Corollary 3.9]{Milne}, therefore $S_3$ is thin.

Now if $\ell \in \pdual(K) \setminus  \cup_{1 \leq i \leq 3} S_i$
then $\ell \cdot X$ is an irreducible effective divisor of
degree $4$, and so there is some quartic extension $L/K$
and a point $P \in X(L)$ such that $\ell$ passes through $P$.
Moreover, $P$ is a non-cusp, and corresponds to a pair
$(E^\prime, u)$ where $E^\prime$ is an elliptic curve over $L$
and $u : E[7] \rightarrow E^\prime[7]$ is a symplectic isomorphism
of $G_L$-modules. 
This supplies an extension $L$ and elliptic curve $E^\prime$
satisfying (ii) and (iv). 
To obtain condition (iii), we need to remove another thin set.
When (iii) fails, $K(E[7]) \cap L$ is either equal to $L$
or is a quadratic extension $M/K$ contained in $K(E[7])$. 
The former is impossible by (b) above. To eliminate the
latter possibility, note that there are at most
finitely many $M \subset K(E[7])$ that are quadratic extensions
of $K$, and for each of these we want our line to avoid the image of the natural
map
\[
X^{(2)}(M) \rightarrow X^{(4)}(K), \qquad D\mapsto D+D^\sigma,
\]
where $\sigma$ denotes conjugation. Of course we have a decomposition
of $X^{(2)}(M)$ similar to \eqref{eqn:x2}, where the 
$M$ replaces $K$, and the elliptic curves $E_i$ are defined over $M$.
By mimicking the above argument with $E_i \times E_i^\sigma$ replacing
$E_i \times E_j$, we deduce the existence of a thin set $S_4$
such that for $\ell \in \pdual$ avoiding $S_1,\dots,S_4$, the conditions
(ii), (iii) and (iv) hold.

We now look at condition (v).
As before, we use the fact that the maximal subgroups of $\GL_2(\F_3)$ are $\SL_2(\F_3)$,
the Borel subgroups, and the normalizers of non-split Cartan subgroups.
If the image of  $\overline{\rho}_{E^\prime, 3}$ does not contain
$\SL_2(\F_3)$ then $E^\prime$
gives 
rise to an $L$-point on either of the geometrically irreducible
curves $Y_1$ or $Y_2$, which are the normalizations of $X \times_{X(1)} X(\mathrm{b}3)$ or $X \times_{X(1)} X({\mathrm{ns}}3)$ respectively.
In particular, $\ell \cdot X$ belongs to the image of 
$Y^{(4)}(K) \rightarrow X^{(4)}(K)$ where $Y=Y_1$ or $Y_2$.
By Lemma~\ref{lem:cover} below, such $\ell$ 
respectively belong 
to thin sets $S_5$, $S_6$. Choosing $\ell$ outside $S_1,\dots,S_6$
yields an $L$ and $E^\prime$ satisfying (ii)--(v).

Finally we want to show that we can choose $\ell$ so that $L$
is totally real. Let $\sigma_1,\dots,\sigma_n$ be the 
real places of $K$. Write $X_{\sigma_i}/\R$ for $X \times_K K_{\sigma_i}$.
Then $X_{\sigma_i}$ is a real twist of the Klein quartic \eqref{eqn:klein}.
By Lemma~\ref{lem:klein} below, the $X_{\sigma_i}$ are in fact
isomorphic to $X(7)$ over $\R$.
The line $2y+z=0$ meets the model given for $X(7)$ in \eqref{eqn:klein} in four distinct real points. Thus
there are non-empty open subsets $U_i \subset \pdual(K_{\sigma_i})$
such that any $\ell \in U_i$ meets $X_{\sigma_i}$ in four
real points. Let
\[
W=\pdual(K) \cap \prod_{i=1}^n U_i.
\]
Now if $\ell \in W \setminus \cup S_i$
then the field $L$ is totally real. All that remains
is to show that $W$ is not thin, and so the
difference $W \setminus \cup S_i$ is non-empty. Note that as $\prod_{i=1}^n \pdual(\R)$ is compact, it is covered by finitely many translates of $\prod_{i=1}^n U_i$ under $\prod_{i=1}^n \PGL_3(\R)$. As $\PGL_3(K)$ is dense in $\prod_{i=1}^n\PGL_3(\R)$ by weak approximation, we can assume the translates are given by elements in $\PGL_3(K)$. Thus $\pdual(K)$ is the union
of finitely many translates of $W$ by the action
of $\PGL_3(K)$. By Hilbert's irreducibility, $W$ is not thin and so $W \setminus \cup S_i$
is non-empty.
\end{proof}

\begin{lem}\label{lem:cover}
Let $X \subset \PP^2$ be a smooth quartic curve over
a number field $K$
and 
 $\mu : Y \rightarrow X$ be $K$-morphism of degree at
least $2$.
Write $\mu$ also for the induced map $Y^{(4)} \rightarrow X^{(4)}$.
Let $\phi: \pdual \hookrightarrow X^{(4)}$ be the map
sending a line $\ell$ to the zero-cycle $\ell \cdot X$.
Then the induced map $\mu^{-1} (\phi(\pdual)) \rightarrow \phi(\pdual)$
does not have a rational section.
\end{lem}
\begin{proof}
By Riemann--Hurwitz, the genus of $Y$ is at least $5$.
We argue by contradiction. Suppose the map has a rational section. Then we have
a commutative triangle,
\[
\xymatrix{
& Y^{(4)} \ar_{\mu}[d]\\
\pdual \ar@{-->}^{\psi}[ru]
\ar@{^{(}->}_{\phi}[r]
& X^{(4)} \\
}
\]
where $\psi$ is a $K$-rational map with $\pdual$ birational to its image.
Choose  $D \in Y^{(4)}(K)$ belonging to the 
image of $\pdual(K)$ under  $\psi$.
Note that $\psi(\pdual)$ is a rational surface
and so its image in the Jacobian $J_Y$ of $Y$ under a suitable
Abel-Jacobi map $\psi(\pdual) \subset Y^{(4)} \rightarrow J_Y$
must be a single point. Thus the
closed points of $\psi(\pdual)$ (interpreted as divisors on $Y$)
are linearly equivalent to $D$. It follows that $\PP(L(D))$ has dimension
at least $2$, and so $l(D) \ge 3$. 
An easy application of Riemann--Roch shows that $D$
is a special divisor (i.e.\ $l(K_Y-D)>0$
where $K_Y$ is the canonical divisor on $Y$).
Clifford's
Theorem \cite[Theorem IV.5.4]{Hartshorne} states that
\begin{equation}\label{eqn:clifford}
l(D)-1 \le \frac{\deg(D)}{2},
\end{equation}
and that if equality holds then $D$ is either canonical, or $Y$ 
is hyperelliptic.
Since $\deg(D)=4$, we see that $l(D)=3$
and so the Clifford inequality~\eqref{eqn:clifford} is in fact an equality.
Thus $D$ is either canonical, or $Y$ is hyperelliptic.
In the first case the genus of $Y$ is $(\deg(D)+2)/2=3$
giving a contradiction. In the second case, $Y$ is a hyperelliptic curve dominating $X$, therefore $X$ is also hyperelliptic by \cite[Proposition 1]{HS}, and we also have a contradiction.
\end{proof}

The following lemma appears to be implicitly assumed
in Manoharmayum's argument \cite[page 710]{Mano1}, and
is perhaps a known fact about the Klein quartic, although
we have been unable to find it in the literature, and therefore give a proof.
\begin{lem}\label{lem:klein}
The only real
twist of the Klein quartic is itself.
\end{lem}
\begin{proof}
Write $X=X(7)$ for the Klein quartic \eqref{eqn:klein}.
Its real twists are classified by $\mathrm{H}^1( \Gal(\C/\R) , \Aut(X))$,
where $\Aut(X)$ is the automorphism group of $X/\C$.
The group $\Aut(X)$ is the unique simple group of order $168$.
Generators for $\Aut(X)$ as a subgroup
of $\PGL_3(\C)$ were given by Klein, and are reproduced by Elkies
in \cite[page 54]{Elkies}. We used these to explicitly write down the semidirect product
\[
\Gal(\C/\R) \ltimes \Aut(X).
\]
There are $28$ complements of $\Aut(X)$ in this semidirect product
and we verified that they form a single conjugacy class, therefore
$\mathrm{H}^1(\Gal(\C/\R) , \Aut(X) ) =0$.
\end{proof}

\section{On Modularity over Totally Real Fields: Proof of Theorem~\ref{thm:Karb}}\label{sec:Karbproof}
Let $K$ be a totally real
field and $E$ an elliptic curve over $K$. By Theorems~\ref{thm:35} 
and~\ref{thm:7}
we know that $E$ is modular
except possibly if $\overline{\rho}_{E,p}(G_{K(\zeta_p)})$ is 
simultaneously absolutely reducible
for $p=3$, $5$, $7$. It follows from Proposition~\ref{prop:large} 
that such $E$ gives rise to a $K$-point on one of the $27$ modular curves
in \eqref{eqn:27mod}.
These all have genera $>1$, and so by Faltings' Theorem~\cite{Faltings}
at most finitely many $K$-points. Theorem~\ref{thm:Karb} follows
immediately.

\bigskip

\noindent \textbf{Remark.} 
By studying the gonalities of various modular curves,
the second-named author \cite{Bao} proves much stronger finiteness statements
than Theorem~\ref{thm:Karb}.

\part{A More Precise Study of Images}

\section{The Image $\overline{\rho}_{E,p}(G_{K(\zeta_p)})$ Again}
\label{sec:imagesagain}

The following proposition gives a refinement of Proposition~\ref{prop:large}
for $p=3$, $5$, $7$. 
Part (a) of the proposition (corresponding to $p=3$) 
is in fact proved by Rubin in \cite[Proposition 6]{Rubin}.
The group denoted by $H$ in part (b) appears (implicitly) 
in a similar setting in \cite[Proof of Lemma 7.2.3]{CDT}.
\begin{prop}\label{prop:large2}
Let $K$ be a totally real number field,
and let $E$ be an elliptic curve over $K$. 
Let $p$ be a rational prime and suppose 
\[
K \cap \Q(\zeta_p)=\Q.
\]
Write 
$\overline{\rho}=\overline{\rho}_{E,p}$. 
Suppose $\overline{\rho}$ is irreducible but 
$\overline{\rho}(G_{K(\zeta_p)})$ is absolutely reducible.
\begin{enumerate}
\item[(a)] If $p=3$ then $\overline{\rho}(G_K)$
is conjugate to $C_{\mathrm{s}}^+(3)$ (the normalizer of a split Cartan subgroup of $\GL_2(\F_3)$);
this has order $8$, and its image in $\PGL_2(\F_3)$ is isomorphic to $(\Z/2\Z)^2$.
\item[(b)] If $p=5$ then $\overline{\rho}(G_K)$
is conjugate in $\GL_2(\F_5)$ to the subgroup
\[
H=\Biggl\langle
\begin{pmatrix}
0 & 1 \\ 2 & 0
\end{pmatrix}, 
\begin{pmatrix}
1 & 0\\
0 & 4
\end{pmatrix}
\Biggr\rangle
\]
of order $16$.
The subgroup $H$ is contained as a subgroup of index $2$
in the normalizer of a split Cartan subgroup,
and as a subgroup of index $3$ in the normalizer of a
non-split Cartan subgroup. The image of $H$ in $\PGL_2(\F_5)$
is isomorphic to $(\Z/2\Z)^2$.
\item[(c)] If $p=7$ then $\overline{\rho}(G_K)$
is conjugate in $\GL_2(\F_7)$ to a subgroup of either
$H_1$ or $H_2$ where
\[
H_1=
\Biggl\langle
\begin{pmatrix}
3 & 0\\
0 & 5
\end{pmatrix},
\begin{pmatrix}
0 & 2\\
2 & 0
\end{pmatrix}
\Biggr\rangle,
\qquad
H_2=
\Biggl\langle
\begin{pmatrix}
0 & 5\\
3 & 0
\end{pmatrix},
\begin{pmatrix}
5 & 0\\
3 & 2
\end{pmatrix}
\Biggr\rangle.
\]
The group $H_1$ has order $36$ and is contained as a subgroup of index $2$ in
the normalizer of a split Cartan subgroup. The group $H_2$
has order $48$ and
is contained as a subgroup of index $2$ in the normalizer
of a non-split Cartan subgroup.
The images of $H_1$ and $H_2$ in $\PGL_2(\F_7)$
are respectively isomorphic to $D_3 \cong S_3$ and $D_4$.
\end{enumerate}
\end{prop}
\begin{proof}
Write $G=\overline{\rho}(G_K)$. We observe that $G$ satisfies
the following:
\begin{itemize}
\item[(i)] $G$ is irreducible, and $G \cap \SL_2(\F_p)$
is absolutely reducible;
\item[(ii)] $G$ odd (in other words, 
there is some element $c \in G$ with eigenvalues $1$ and $-1$); 
\item[(iii)] $\det(G)=\F_p^*$.
\end{itemize}
Indeed, (i) follows by assumption and Proposition~\ref{prop:large},
(ii) is a consequence of $K$ having a real embedding, 
and (iii) follows from $K \cap \Q(\zeta_p)=\Q$.
One proof of the proposition is to enumerate all subgroups
$G$ of $\GL_2(\F_p)$ with $p=3$, $5$, $7$,
and check that those satisfying (i), (ii), (iii),
also satisfy (a), (b), (c).
We can simplify the enumeration by using Lemma~\ref{lem:exceptional}
and its proof, from which we know that $G \cap \GL_2^+(\F_p)$
is contained in a Cartan subgroup $C_*(p)$ and $G$
is contained in its normalizer $C_*^+(p)$.
If $p=5$ then the element $c \in G$ that represents complex
conjugation belongs to $\GL_2^+(\F_p)$ and therefore to
$G \cap \GL_2^+(\F_p)$. Now an element of $C_{\mathrm{ns}}(5)$
has eigenvalues $\alpha$, $\alpha^p$ for some $\alpha \in \F_{25}^*$.
This shows that the non-split Cartan case does not arise for $p=5$.
\end{proof}
\noindent{\bf Remark}. It follows from the proof that $\overline{\rho}(G_{L})$
is contained in a Cartan subgroup, where $L=K(\sqrt{(-1)^{(p-1)/2} p})$.

\section{Modular Curves Again}
\label{sec:modcurvesagain}
Let $H_1$, $H_2$ be the
two subgroups of $\GL_2(\F_7)$ specified
in part (c) of Proposition~\ref{prop:large2}.
We shall denote the corresponding modular curves $X(H_1)$ and $X(H_2)$
by $X(\mathrm{d}7)$ and $X(\mathrm{e}7)$ respectively.
We extend our earlier conventions regarding the naming
of modular curves in an obvious way. Thus, for a number field $K$,
non-cuspidal $K$-points on $X(\mathrm{s}3,\mathrm{b}5,\mathrm{d}7)$
give elliptic curves $E/K$ such that the image of $\overline{\rho}_{E,3}$
is contained in a conjugate of the normalizer of a split Cartan subgroup of $\GL_2(\F_3)$, the image of
$\overline{\rho}_{E,5}$ is contained in a conjugate of the Borel subgroup of $\GL_2(\F_5)$,
and the image of $\overline{\rho}_{E,7}$ is contained in a conjugate of the group $H_1$.

\begin{cor}\label{cor:large2}
Let $K$ be a real quadratic field,
and let $E$ be an elliptic curve defined over $K$.
\begin{enumerate}
\item[(a)] If $\overline{\rho}_{E,3} (G_{K(\zeta_3)})$
is absolutely reducible then $E$ gives rise to a non-cuspidal
$K$-point on either $X(\mathrm{b}3)$ or $X(\mathrm{s}3)$.
\item[(b)] If $K \ne \Q(\sqrt{5})$ 
and $\overline{\rho}_{E,5} (G_{K(\zeta_5)})$
is absolutely reducible then $E$ gives rise to a non-cuspidal
$K$-point on either $X(\mathrm{b}5)$ or $X(\mathrm{s}5)$.
\item[(c)] If $\overline{\rho}_{E,7} (G_{K(\zeta_7)})$
is absolutely reducible then $E$ gives rise to a non-cuspidal
$K$-point on either $X(\mathrm{b}7)$ or $X(\mathrm{d}7)$ or $X(\mathrm{e}7)$.
\end{enumerate}
\end{cor}
\begin{proof}
The field $K$ is real quadratic, so for $p=3$, $7$, we have
$K \cap \Q(\zeta_p)=\Q$. For $p=5$ this holds once we have
imposed the condition $K \ne \Q(\sqrt{5})$. The corollary
follows from Proposition~\ref{prop:large2}.
\end{proof}
\noindent{\bf Remarks}.
\begin{enumerate}
\item[(i)] In (b) we could have replaced $X(\mathrm{s}5)$
with the curve $X(H)$, where $H$ is as given in 
the statement of Proposition~\ref{prop:large2}.
In fact $X(H)$ is a double cover of $X(\mathrm{s}5)$,
so this would give a more precise statement.
However, for our later purpose of proving
modularity of elliptic curves
over real quadratic fields, we have found
it sufficient and simpler to work with
$X(\mathrm{s}5)$. 
\item[(ii)] Parts (a), (b) improve over
Proposition~\ref{prop:large} by eliminating
$X(\mathrm{ns}3)$ and $X(\mathrm{ns}5)$.
This substantially reduces the amount of work
we have to do later in enumerating quadratic
points on modular curves.
\item[(iii)] In (b) we really need the condition $K \ne \Q(\sqrt{5})$.
To see this let $K=\Q(\sqrt{5})$. Let $G=\overline{\rho}(G_K)$
where $\overline{\rho}=\overline{\rho}_{E,5}$ and $E/K$
is an elliptic curve. Then 
\[
\det(G)=\chi_5\left(\Gal(\Q(\zeta_5)/\Q(\sqrt{5})) \right) =
\{ \overline{1}, \overline{4} \} \subset \F_5^*.
\]
Now, in the proof of Proposition~\ref{prop:large2},
if we replace the condition
$\det(G)=\F_5^*$ with $\det(G)=\{\overline{1},\overline{4}\}$
we obtain three subgroups of $\GL_2(\F_5)$ of orders
$6$, $8$ and $12$. The subgroup of order $8$ is contained
in the normalizer of a split Cartan subgroup, but not 
in the normalizer of any non-split Cartan subgroup.
The subgroups of orders $6$ and $12$ are 
contained in the normalizer of a non-split Cartan
subgroup but not in the normalizer
of any split Cartan subgroup.
\end{enumerate}

\section{Theorem~\ref{thm:modquadratic} follows from Lemmas~\ref{lem:Qsqrt5}
and~\ref{lem:7mod}}
\label{sec:modcurvesquad}

Most of the rest of this paper is devoted to the proof
of Lemma~\ref{lem:7mod}, which states that non-cuspidal
real quadratic points on the modular curves in \eqref{eqn:first3}
and \eqref{eqn:last4} are modular.
Let $X$ be any of the modular curves in \eqref{eqn:first3},~\eqref{eqn:last4}.
Recall that $X$ is a moduli space of elliptic curves with level
structures. A non-cuspidal point on $X$ gives rise to a well-defined
$j$-invariant, and when $j\neq 0, 1728$ this determines the underlying elliptic
curve up to quadratic twists. When we say that a point on $X$ is modular, we 
mean that one (and hence any) elliptic curve supporting it is modular.
Lemma~\ref{lem:7mod}
is in fact a summary
of results established separately in the following lemmas:
Lemma~\ref{lem:b5b7}, Lemma~\ref{lem:b3s5},
Lemma~\ref{lem:s3s5} and Lemma~\ref{lem:b3b5d7}.

Lemma~\ref{lem:Qsqrt5} asserts modularity
of elliptic curves over $\Q(\sqrt{5})$. 
Its proof is relegated to Section~\ref{sec:Qsqrt5}. 
We now explain how Lemmas~\ref{lem:Qsqrt5} and \ref{lem:7mod}
are enough to imply the modularity of all elliptic curves over real
quadratic fields.
\begin{proof}[Proof of Theorem~\ref{thm:modquadratic}]
Let $E$ be an elliptic curve over a real quadratic field $K$.
We would like to show that $E$ is modular. By Lemma~\ref{lem:Qsqrt5}
we may suppose that $K \ne \Q(\sqrt{5})$.
We know from Theorems~\ref{thm:35} and~\ref{thm:7}
that if $\overline{\rho}_{E,p}(G_{K(\zeta_p)})$ is absolutely
irreducible with $p=3$ or $5$ or $7$ then $E$ is modular.
Thus we may restrict
our attention to elliptic curves $E$ such that
$\overline{\rho}_{E,p}(G_{K(\zeta_p)})$ is 
simultaneously absolutely reducible
for $p=3$, $5$ and $7$. 
Now appealing to Corollary~\ref{cor:large2} we see that
$E$ gives rise to a $K$-point on one of the twelve modular curves
\begin{gather*}
X(\mathrm{b}3,\mathrm{b}5,\mathrm{b}7),
\qquad
X(\mathrm{b}3,\mathrm{b}5,\mathrm{d}7),
\qquad
X(\mathrm{b}3,\mathrm{b}5,\mathrm{e}7),\\
X(\mathrm{b}3,\mathrm{s}5,\mathrm{b}7),
\qquad
X(\mathrm{b}3,\mathrm{s}5,\mathrm{d}7),
\qquad
X(\mathrm{b}3,\mathrm{s}5,\mathrm{e}7),\\
X(\mathrm{s}3,\mathrm{b}5,\mathrm{b}7),
\qquad
X(\mathrm{s}3,\mathrm{b}5,\mathrm{d}7),
\qquad
X(\mathrm{s}3,\mathrm{b}5,\mathrm{e}7),\\
X(\mathrm{s}3,\mathrm{s}5,\mathrm{b}7),
\qquad
X(\mathrm{s}3,\mathrm{s}5,\mathrm{d}7),
\qquad
X(\mathrm{s}3,\mathrm{s}5,\mathrm{e}7).\\
\end{gather*}
It is clear that every one of these twelve covers one of the seven 
modular curves in Lemma~\ref{lem:7mod}, which completes the proof.
\end{proof}

\bigskip

The proof of Lemma~\ref{lem:7mod} will involve
the determination of all non-cuspidal real quadratic points
on some of the seven modular curves. If these points
correspond to CM elliptic curves or $\Q$-curves then
they are modular (see the remark at the beginning of
Section~\ref{sec:b5b7} for more on the latter).
The following lemma will sometimes be useful in proving
the modularity of the remaining points.
\begin{lem}\label{lem:bigimage}
Let $E$ be an elliptic curve over a totally real field $K$,
with $j$-invariant $j$. Let $p=3$ or $5$ or $7$, and suppose
that the $p$-th division polynomial of $E$ is irreducible.
Let $\fq$ be a prime of $K$ satisfying $\upsilon_\fq(j)<0$
and $p \nmid \upsilon_\fq(j)$. Then $E$ is modular.
\end{lem}
\begin{proof}
We shall show that $\overline{\rho}_{E,p}(G_K)$ contains $\SL_2(\F_p)$,
which is enough to imply modularity by Theorems~\ref{thm:35}
and~\ref{thm:7}.

By the theory of the Tate curve \cite[Proposition 6.1]{SilvermanII},
the two conditions $\upsilon_\fq(j)<0$
and $p \nmid \upsilon_{\fq}(j)$ together imply 
that $p$ divides the order of $\overline{\rho}_{E,p}(G_K)$.
It follows from the classification of maximal subgroups
of $\GL_2(\F_p)$ (see for example \cite[Lemma 2]{SwD})
that $\overline{\rho}_{E,p}(G_K)$ is either contained
in a Borel subgroup or contains $\SL_2(\F_p)$.
Since the
$p$-th division polynomial is irreducible, 
$\overline{\rho}_{E,p}(G_K)$ is not contained
in a Borel subgroup. 
\end{proof}

\noindent \textbf{Remarks}.

\medskip

\noindent \textbf{(I)}
See \cite{Bao}
for an alternative study of modular curves and quadratic points on them related
to elliptic curves with small residual image at $3$, $5$, and $7$.

\medskip

\noindent \textbf{(II)} It is natural
to ask whether the quadratic points on the seven modular curves 
in Lemma~\ref{lem:7mod} can determined by
the methods of Mazur \cite{MazurEisen}, Kamienny \cite{Kamienny} and 
Merel \cite{Merel}. 
Let $X$ be any of these modular curves.
To apply such methods to our problem we appear to require:
\begin{enumerate}
\item[(a)] an abelian variety $\cA/\Q$ having Mordell--Weil rank $0$,
that is a quotient
of the Jacobian $J$ of $X$; 
\item[(b)] a prime $p$ and a morphism $X^{(2)} \rightarrow \cA$ (defined
over $\Q$) that is
injective on some small $p$-adic neighbourhood $U \subset X^{(2)} (\Q_p)$ (where $X^{(2)}$ is the symmetric square 
of $X$);
\item[(c)] some argument that shows that $X^{(2)}(\Q) \subset U$. 
\end{enumerate}
Usually $U$ is the $p$-adic 
residue class of $\{ \infty, \infty\}$, and property (b) 
is a consequence of $X^{(2)} \rightarrow \cA$ being a formal immersion
at $\{\infty,\infty\}$ above $p$. For the modular curves that we consider,
we have been unable to demonstrate (c) without actually determining the
quadratic points~\footnote{
Indeed, the absence of (c) is the reason why 
we do not yet have a parametrization
of quadratic points on the family $X(\mathrm{b}N)$, even though
(a) and (b) are known \cite{Kamienny} for prime $N>71$, where $\cA$
is the Eisenstein quotient of $J_0(N)$, and 
$p$ is any prime $\ne 2$, $3$, $5$.}.
 In fact, as we will see in due course,
$X(\mathrm{b}5,\mathrm{b}7)$
and $X(\mathrm{s}3,\mathrm{s}5)$ have infinitely many quadratic
points, so (b) and (c) cannot simultaneously hold. Likewise,
$X(\mathrm{b}3,\mathrm{s}5)$ has many non-cuspidal
quadratic points, suggesting that these methods are inapplicable.

\part{Real Quadratic Points on Certain Modular Curves}

\section{Real Quadratic Points on $X(\mathrm{b}5,\mathrm{b}7)$}\label{sec:b5b7}
In this section we parametrize the non-cuspidal real quadratic points
on $X(\mathrm{b}5,\mathrm{b}7)$ (commonly known as $X_0(35)$), 
showing that they are all
modular. In fact, we shall show that all except two correspond to $\Q$-curves.
By a \text{$\Q$-curve}, we mean an elliptic  curve defined over $\overline{\Q}$
that is isogenous to all its Galois conjugates. 
It is known that
$\Q$-curves are modular. Indeed,
Ribet  \cite[Corollary 6.2]{Ribet} had shown that modularity of $\Q$-curves
is a consequence of Serre's modularity conjecture,
which is now a theorem of Khare and Wintenberger \cite{KW1}, \cite{KW2}.
\begin{lem}\label{lem:b5b7}
Let $X=X(\mathrm{b}5,\mathrm{b}7)$. The non-cuspidal real
quadratic points on $X$ are modular.
\end{lem}
\begin{proof}
The curve $X$ is hyperelliptic of genus $3$. 
Galbraith \cite[Section 4.4]{Galbraith} derives the following model for $X$,
\begin{equation}\label{eqn:35}
X \; : \; y^2=(x^2 + x - 1)(x^6 - 5 x^5 - 9 x^3 - 5 x - 1).
\end{equation}
Write $J$ for $J_0(35)$---the Jacobian of $X$. Using $2$-descent
\cite{Stoll} we find that $J(\Q)$ has Mordell-Weil rank $0$.
As $J$ has good reduction at $3$, the map $J(\Q) \rightarrow J(\F_3)$
is injective. The group generated by differences of the four obvious
rational points $(0,\pm 1)$, $\infty_{\pm}$ (where the latter
are the points at infinity) surjects onto
$J(\F_3) \cong \Z/24\Z \times \Z/2\Z$, and is therefore equal
to $J(\Q)$. 
We find
\[
J(\Q)= \frac{\Z}{24 \Z} \cdot \left[ \infty_{-}-\infty_{+} \right]
+\frac{\Z}{2 \Z} \cdot \left[3(0,-1) -3 \infty_{+} \right].
\]

Let $P$ be a quadratic point of $X$, and let $P^\sigma$
denote its conjugate. Then $P+P^\sigma-\infty_{+}-\infty_{-}$
is a rational divisor of degree $0$ and hence defines
an element of $J(\Q)$. It follows that $P+P^\sigma$
is linearly equivalent to $T+\infty_{+}+\infty_{-}$
where $T$ is a divisor of degree $0$ representing
one of the $48$ elements of $J(\Q)$. 
As $X$ is hyperelliptic, for $[T] \ne 0$,
the divisor $T+\infty_{+}+\infty_{-}$ is linearly equivalent
to at most one effective divisor of degree $2$. Enumerating 
these $47$ possibilities for $[T] \ne 0$, we find that all
the effective degree $2$ divisors they yield are sums
of two rational points, 
except for $\left(\frac{-1\pm \sqrt{5}}{2},0\right)$. 
We would like to show that this pair of points is modular.
By Ogg \cite[Corollary 2]{ogg},
the hyperelliptic involution on $X=X(\mathrm{b}5,\mathrm{b}7)$
is the Fricke involution $w_{35}$. 
The points $\left(\frac{-1\pm \sqrt{5}}{2},0\right)$ are fixed by
$w_{35}$ and so correspond to CM elliptic curves,
which are therefore modular.

We now turn to the case $[T]=0$, and so $P+P^\sigma$
is linearly equivalent to $\infty_{+}+\infty_{-}$.
Thus
$P=(x,\pm \sqrt{f(x)})$ with
$x \in \Q\backslash \{0\}$, where $f$ is the polynomial on
the right-hand side of \eqref{eqn:35}.
Suppose now that
$P$ is a quadratic point of the form $(x,\sqrt{f(x)})$ with
$x \in \Q\backslash \{0\}$. The Galois conjugate $P^\sigma$ of $P$
is related to $P$ by the hyperelliptic involution. 
Therefore the elliptic
curve $E$ represented by $P \in X$ is related to
its Galois conjugate $E^\sigma$ by a $35$-isogeny, and hence is a
$\Q$-curve.
\end{proof}

\section{Real Quadratic Points on $X(\mathrm{b}3,\mathrm{s}5)$}\label{sec:b3s5}

\begin{lem}\label{lem:b3s5}
Let $X=X(\mathrm{b}3,\mathrm{s}5)$.
All non-cuspidal real quadratic points on $X$ 
are modular.
\end{lem}
\begin{proof}
The curve $X$ is studied briefly in \cite[page 554--556]{CDT}
where it is shown that 
$X \cong X_0(75)/w_{25}$. 
That paper does not give a model for $X$ but
mentions a computation by Elkies where he finds
$X$ to have genus $3$ and determines its rational points.
It turns out there is precisely one non-cuspidal rational
point and that this has $j$-invariant $0$.

In order to write down a model for $X$, and to compute
the Mordell--Weil group 
of its Jacobian, we shall make use of 
Stein's {\tt Magma} implementation of modular symbols
algorithms \cite{Cre}, \cite{Stein}. 
We find that an echelon basis for the 
$+1$-eigenspace of the Atkin--Lehner operator $w_{25}$ acting on
$S_2(\Gamma_0(75))$ is given by
\begin{align*}
f_1 &=    q - q^4 - 2 q^5 + q^6 + q^9 + 2 q^{10} - 4 q^{11} + \cdots\\
f_2 &=    q^2 - q^4 + q^6 - q^7 - q^8 - 2 q^{11} + \cdots\\
f_3 &=    q^3 + q^4 + 2 q^5 - q^6 + q^7 - 2 q^8 - 2 q^{10} + 2 q^{11} + \cdots 
\end{align*}
Thus $X$
has genus $3$; a basis for the holomorphic differentials is
given by $f_i dq/q$. Following Galbraith \cite[Chapter 3]{Galbraith},
we compute the image of this curve under the canonical map,
and find that it is a plane quartic curve. 
It follows that the canonical map is indeed an embedding
of $X$ 
in $\PP^2$. The equation given by the canonical embedding is
\begin{multline}\label{eqn:Xmodel}
X \; :\;
 3 x^3 z - 3 x^2 y^2 + 5 x^2 z^2 - 3 x y^3 
- 19 x y^2 z - x y z^2 \\ + 3 x z^3 + 2 y^4 +
    7 y^3 z - 7 y^2 z^2 - 3 y z^3=0.
\end{multline}
We shall also need the map $j:X \rightarrow X(1)$.
The details of finding this map explicitly are similar
to those in a recent paper of Banwait and Cremona \cite[Section 5]{BC},
and so we only give a sketch.
The isomorphism $X_0(75)/w_{25} \rightarrow X$ 
pulled back to the upper half-plane $\mathfrak{H}$ is given
by $\tau \mapsto 5\tau$. Thus the $q$-expansion
of $j$ on $X$ is given by
\[
j(q^5)=\frac{1}{q^5}+
744 + 196884 q^5 + 21493760 q^{10} + 864299970 q^{15} + 
    20245856256 q^{20} +  \cdots.
\]
We can express $j=F(x,y,z)/G(x,y,z)$ where $F$, $G$
are homogeneous forms of the same degree with coefficients in $\Z$. 
We seek such $F$, $G$ so that the identity of modular forms
\begin{equation}\label{eqn:cuspexp}
F(f_1,f_2,f_3)=j(q^5) G(f_1,f_2,f_3)
\end{equation}
holds. To find suitable $F$, $G$, we run through the degrees $d=1,2,\dots$
and express $F$, $G$ as linear combinations of the monomials of degree $d$
with unknown coefficients. By comparing the coefficients of the $q$-expansions
in \eqref{eqn:cuspexp}
we obtain simultaneous equations for the unknown coefficients of $F$, $G$.
Solving these we obtain the first non-zero $F$, $G$ with $d=15$. 
For this particular $F$, $G$ we checked the equality
of coefficients on both sides of \eqref{eqn:cuspexp} beyond
the Sturm bound \cite[Theorem 9.18]{Stein}. 
This is sufficient to show that $j=F(x,y,z)/G(x,y,z)$.
The expressions for $F$ and $G$ are complicated
and can be found in the accompanying {\tt Magma} script (see 
page~\pageref{magma}).

\medskip

Write $J$ for the Jacobian of $X$. Computing the 
eigenbasis for the space spanned by $f_1$, $f_2$, $f_3$
yields three rational eigenforms corresponding to the
elliptic curves $E_1$, $E_2$, $E_3$ with Cremona  
labels $15A1$, $75B1$, $75C1$. It follows that
$J$ is isogenous to $E_1 \times E_2 \times E_3$. 

We use this to compute the rank of $J(\Q)$
and of $J(\Q(\sqrt{5}))$ (we will see in due course
why the latter rank is convenient to have).
The three curves $E_i$ have rank $0$ over $\Q$
and so do their quadratic twists by by $5$.
Thus $E_i$ have rank $0$ over $\Q$
and over $\Q(\sqrt{5})$.
It follows that $J$ has rank $0$ over those two fields. 
Thus $J(\Q)$
is equal to its torsion subgroup. 
Let $L=\Q(\sqrt{5})$. 
We shall denote the torsion subgroups of $J(\Q)$ and
of $J(L)$ by $J(\Q)_{\mathrm{tors}}$ and
$J(L)_{\mathrm{tors}}$. 
To compute $J(\Q)_{\mathrm{tors}}$ we shall in fact
compute $J(L)_{\mathrm{tors}}$ and recover
$J(\Q)_{\mathrm{tors}}$ using the inclusion
$J(\Q)_{\mathrm{tors}} \subseteq J(L)_\mathrm{tors}$.
A short search for points on $X$ defined over $L$
reveals the following points:
\begin{eqnarray*}
 P_0=(0:0:1),\qquad
 P_1=(1:0:0), \qquad
 P_2=( 1: -1: 2 ),\\
 P_3=\left( -2: \frac{- 1+\sqrt{5}}{2}: 1 \right),\quad
   P_4= \left( \frac{1-\sqrt{5}}{2}: -3+\sqrt{5} : 1 \right),\quad
    P_5=\left( \frac{1+\sqrt{5}}{2}: 2 + \sqrt{5}: 1 \right).
\end{eqnarray*}
As $J(L)$ is torsion, we see that $[P_i-P_0]$
belong to $J(L)_{\mathrm{tors}}$. 
Write $H$ for the subgroup of $J(L)_{\mathrm{tors}}$
generated by the differences $[P_i-P_0]$, $[P_i^\sigma-P_0]$
where $\sigma$ denotes conjugation in $L$.
Write $\OO_L$ for the ring
of integers of $L$ and let $\fp=19 \OO_L+ (9+\sqrt{5})\OO_L$.
Then $\fp$ is a prime ideal of norm $19$. Thus we have a natural
injection
\[
\pi_\fp: J(L)_{\mathrm{tors}} \hookrightarrow J(\F_{19})
\]
induced by $\fp$. We find that
\[
J(\F_{19}) \cong 
\Z/2\Z \times \Z/2\Z \times \Z/8\Z \times \Z/200\Z,
\qquad \pi_\fp(H) \cong 
\Z/2\Z \times \Z/2\Z \times \Z/8\Z \times \Z/40\Z.
\]
It follows that either $J(L)_{\mathrm{tors}}=H$ or 
$H$ has index $5$ in $J(L)_{\mathrm{tors}}$. In the latter case,
the order of $J(L)_{\mathrm{tors}}$ is divisible by $25$.

Let $\fq=7\OO_L$.
This is a prime ideal of norm $49$, and so we have
another injection
\[
\pi_\fq : J(L)_{\mathrm{tors}} \hookrightarrow J(\F_{49})
\cong (\Z/8\Z)^3 \times \Z/440\Z.
\] 
As the order of $J(\F_{49})$ is not divisible by $25$,
the order of $J(L)_{\mathrm{tors}}$ is not divisible
by $25$. It follows that $J(L)_{\mathrm{tors}}=H$.
Taking $\Gal(L/\Q)$-invariants, we find that
\[
J(\Q)=J(\Q)_{\mathrm{tors}}
=(\Z/40\Z)\cdot [P_2-P_0] + (\Z/2\Z) \cdot [P_3+P_3^\sigma + 3 P_0-5 P_2].
\]

We now enumerate the real quadratic points of $X$.
As $X$ is a plane quartic (and therefore not hyperelliptic), the map
\[
X^{(2)}(\Q) \rightarrow J(\Q), \qquad D \mapsto [D-2P_0]
\]
is injective. 
Note that any element $D \in X^{(2)}(\Q)$
is linearly equivalent to $D^\prime+2P_0$, where $[D^\prime] \in J(\Q)$. 
For each of the $80$ elements $[D^\prime]$ of $J(\Q)$ we computed the Riemann--Roch
space of $D^\prime+2 P_0$ and checked if this divisor is linearly
equivalent to a positive divisor of degree $2$ (in which case
the positive divisor is unique). This allowed us to
enumerate $X^{(2)}(\Q)$ and hence all the rational and 
quadratic points on $X$. It turns out that 
\[
X(\Q)=\{ P_0, P_1, P_2 \}.
\]
The points $P_0$ and $P_1$ are cusps; the point $P_2$ has
$j$-invariant $0$. 
This is consistent with Elkies' computation mentioned
in \cite[page 555]{CDT}.
The complex quadratic points are (up to conjugation)
\[
(- 5+\sqrt{-11} : 0 : 6), \qquad
(-6+2 \sqrt{-15} :  1+\sqrt{-15} : 8),
\]
with $j$-invariants
$-32768$ and $0$ respectively.
Table~\ref{table:Xpts} gives the real quadratic points
(up to conjugation),
and their $j$-invariants; the ones with $j$-invariant
$=\infty$ are of course cusps.

From the table we see that there are $12$ non-cuspidal real quadratic points
(up to conjugation). Of these, seven correspond to elliptic curves
with complex multiplication and are therefore modular. For
the remaining five modularity follows from Lemma~\ref{lem:bigimage}
with $p=7$ and $\fq$ some prime above $2$.
\end{proof}

\begin{landscape}
\begin{table}[h] 
\centering
\begin{tabular}{|c | c | c| }
\hline
point & $j$-invariant & discriminant of quadratic \\
& & order of CM \\
\hline\hline
$(-2+\sqrt{2}: \sqrt{2} :2 )$ &
$2417472 -1707264 \sqrt{2} $ &
$-24$ \\
\hline
$(-1 : - 2+\sqrt{3} : 1)$ &
$1728$ &
$-4$ \\
\hline
$(-2 :  - 1+\sqrt{3} : 2)$ &
$76771008 + 44330496 \sqrt{3}$ &
$-36$ \\
\hline
$( 1+\sqrt{5} : - 6-2\sqrt{5} : 2)$ &
$0$ &
$-3$ \\
\hline
$( 1+\sqrt{5} :  4+2\sqrt{5} : 2)$ &
$\infty$ &
not applicable \\
\hline
$(-4 : - 1+\sqrt{5} : 2)$ & 
$\infty$ & 
not applicable \\
\hline
$( - 10+2\sqrt{17} : 3 - \sqrt{17} : 4)$ &
$- 2770550784 -671956992 \sqrt{17} $ &
$-51$ \\
\hline
$(- 3+\sqrt{33} : 6 : 0)$ &
$-32768$ &
$-11$ \\
\hline
$( - 9+ \sqrt{33} :  - 9+\sqrt{33} : 8)$ &
$- 18808030478336 +3274057859072 \sqrt{33}$ & 
$-99$ \\
\hline
$(-8 : - 11+\sqrt{41} : 4)$ &
$( 3674587331712451 -573874126975281 \sqrt{41} )/ 65536$ &
non-CM \\
\hline
$(- 5+\sqrt{41} : - 2+2\sqrt{41} : 8)$ &
$( - 251512901 -39335769 \sqrt{41} )/64 $ &
non-CM\\
\hline
$(- 23+\sqrt{145} : -12 : 24)$ &
$( 3494273 -283495 \sqrt{145} )/64 $ &
non-CM\\
\hline
$(11+\sqrt{145} : - 11-\sqrt{145} : 2)$ &
$( 839560577 -84563815 \sqrt{145} ) /65536$ &
non-CM \\
\hline
$( - 91+\sqrt{185} : 6-2\sqrt{185} : 88)$ &
$ ( 572227476834375 -42070984346475 \sqrt{185}  )/ 65536 $ &
non-CM \\
\hline
\end{tabular}
\caption{The table gives the real quadratic points
on the curve $X(\mathrm{b}3,\mathrm{s}5)$
given by the model in \eqref{eqn:Xmodel}.
It also gives the $j$-invariants of these points,
and, where applicable, the discriminant
of the order of complex multiplication of 
the corresponding elliptic curve.
}
\label{table:Xpts}
\end{table}
\end{landscape}

\section{Real Quadratic Points on $X(\mathrm{s}3,\mathrm{s}5)$}\label{sec:s3s5}

\begin{lem}\label{lem:s3s5}
Non-cuspidal real quadratic points on $X(\mathrm{s}3,\mathrm{s}5)$
are modular.
\end{lem}
\begin{proof}
Let $X=X(\mathrm{s}3,\mathrm{s}5)$.
As in the proof of Lemma~\ref{lem:b3s5}, we use the canonical
map to write down a model for $X \cong X_0(225)/\langle w_9, w_{25} \rangle$.
The subspace of $S_2(\Gamma_0(225))$ invariant under $w_9$ and $w_{25}$
has dimension $4$ and echelon basis
\begin{align*}
f_1 &= q - 4 q^5 - 4 q^6 + 5 q^7 - 4 q^9 + 4 q^{10} - 8 q^{11} + \cdots,\\
f_2 &=  q^2 - 2 q^5 - 4 q^6 + 4 q^7 - q^8 - 2 q^9 + 2 q^{10} - 6 q^{11} + \cdots,\\
f_3 &=    q^3 - 2 q^5 - 2 q^6 + 2 q^7 + q^8 - q^9 + 2 q^{10} - 3 q^{11} + \cdots,\\
f_4 &=    q^4 - 2 q^5 - 2 q^6 + 5 q^7 - 2 q^9 + 2 q^{10} - 4 q^{11} + \cdots 
\end{align*}
Thus $X$ has genus $4$. Instead of using the basis $f_1,\ldots,f_4$,
we shall use the following basis:
\[
g_1= -f_1+2 f_4,
\quad
g_2= f_1 -f_2+2 f_3 -f_4,
\quad
g_3=f_1+f_2-2 f_3 -f_4,
\quad
g_4= f_1 -f_3.
\]
The motivation behind this change of basis is that the 
automorphisms of the model
for $X$ given by the canonical map induced by $g_1,\dots,g_4$
are easier to describe as we will soon see. The canonical map
$X \rightarrow \PP^3$ induced by $g_1,\dots,g_4$ is an embedding
and yields the following model:
\[
X \; : \; \left\{
\begin{array}{ll}
3 x_1^2 - x_2^2 - 2 x_3^2 + 2 x_3  x_4 - 3 x_4^2=0\\
x_2^2 x_3 - x_3^3 + 4 x_3^2 x_4 - 12 x_3 x_4^2 + 12 x_4^3=0.
\end{array}
\right.
\]
Note that $X$ has the following pair of involutions,
\[
\alpha : (x_1,x_2,x_3,x_4) \mapsto (-x_1,x_2,x_3,x_4), 
\quad
\beta : (x_1,x_2,x_3,x_4) \mapsto (x_1,-x_2,x_3,x_4).
\]
Write
$\cA:=\{1,\alpha,\beta,\alpha\beta\}$.
(It is possible to show that $\cA$
is the automorphism group of $X$, but we shall not need this.) 
Write $\mathbb{K}$ for the function field of $X$.
The subfield fixed by $\cA$ is $\Q(w)$ where $w=x_3/x_4$.
Since $\cA\cong \Z/2\Z \times \Z/2\Z$, it follows that
$\mathbb{K}=\Q(w)(\sqrt{h_1},\sqrt{h_2})$ for some 
$h_1$, $h_2 \in \Q(w)$. By inspection, we see that $X$ has the following
more convenient model:
\begin{equation}\label{eqn:s3s5}
X \; : \; \left\{
\begin{array}{ll}
u^2&=w (w-1) (w^2-w+4)\\
v^2 & =w (w^3-4 w^2+12 w-12),
\end{array}
\right.
\end{equation}
where
\[
u=\frac{x_1  x_3}{ x_4^2},
\qquad
v=\frac{x_2 x_3}{x_4^2},
\qquad
w=\frac{x_3}{x_4}.
\]

We are interested in determining the real 
quadratic points of $X$, and for this we shall
consider the following three obvious quotients
\begin{align*}
D_1 \; &: \;
u^2= w (w-1) (w^2-w+4),\\
D_2 \; &: \;
v^2  =w (w^3-4 w^2+12 w-12),\\
D_3 \; &: \;
z^2 = (w-1) (w^2-w+4)(w^3-4 w^2+12 w-12).
\end{align*}
We shall write
\begin{align*}
h_1& =w (w-1) (w^2-w+4), \qquad \qquad
h_2=w (w^3-4 w^2+12 w-12),\\
h_3 &=(w-1) (w^2-w+4)(w^3-4 w^2+12 w-12).
\end{align*}
The quotient maps $\phi_i : X \rightarrow D_i$ are all of degree $2$,
given by
\[
\phi_1(w,u,v)=(w,u),\qquad
\phi_2(w,u,v)=(w,v),\qquad
\phi_3(w,u,v)=(w,uv/w).
\]
The curve $D_2$ is isomorphic to the Weierstrass elliptic
curve with Cremona label $225A1$, with Mordell--Weil rank $1$ over $\Q$ 
(and trivial torsion). 
It follows that $X$ has infinitely many
quadratic points obtained by pulling back the rational
points on $D_2$ via the degree $2$ map $\phi_2$.
 We claim that all non-cuspidal real quadratic
points on $X$ are obtained thus. 

The curve $D_3$ has genus $2$. Write $J$ for its Jacobian.
Using $2$-descent we find that $J(\Q)$ has Mordell--Weil rank $0$,
and by a similar calculation to those in the previous sections we find that
\[
J(\Q)=\frac{\Z}{10\Z} \cdot \left[
\left(
\frac{1+\sqrt{5}}{2}, \frac{-15+5\sqrt{5}}{2}
\right)
+
\left(
\frac{1-\sqrt{5}}{2}, \frac{-15-5\sqrt{5}}{2}
\right)
-\infty_{+} -\infty_{-}
 \right] \, .
\]
Here $\infty_{+}$, $\infty_{-}$ are the two points at infinity on $D_3$.
As in the proof of Lemma~\ref{lem:b5b7} we can determine all
rational and quadratic points on $D_3$. 
In particular,
\[
D_3(\Q)=\{ (1,0), \infty_+, \infty_{-} \}.
\]
The real quadratic points $(w,z)$ on $D_3$
are, up to Galois conjugation and the action of the hyperelliptic involution,
equal to 
$((1+ \sqrt{5})/2, (-15+ 5\sqrt{5})/2)$,
%
or satisfy $w \in \Q$ and $z=\pm \sqrt{h_3(w)} \notin \Q$.

Clearly a real quadratic point on $X$ maps to either 
a rational or a real quadratic point on $D_3$.
We deduce that the only real quadratic points on $X$ 
satisfy 
\begin{enumerate}
\item[(a)] either $(w,u,v)=
( (1-\sqrt{5})/2, (-5-\sqrt{5})/2, \sqrt{5})$
or its conjugates under the action of Galois and the automorphisms
$u \mapsto -u$, $v \mapsto -v$;
\item[(b)] or $w \in \Q$ and $uv/w=\pm \sqrt{h_3(w)} \notin \Q$,
\end{enumerate}
The points  in (a) satisfy $w^2-w-1=0$, and are cusps as
can be deduced from the expression for $j$-map on $X$ given below
in \eqref{eqn:js3s5}.

We are interested in proving the modularity of the
non-cuspidal real points $P=(w,u,v)$ satisfying (b). Write $d=h_3(w)$.
Let $K=\Q(\sqrt{d})$, which is the field of the definition
of the point $P$. If $u$ is rational, then $\phi_1(P)=(w,u) \in D_1(\Q)$.
However,
$D_1$ is isomorphic to the Weierstrass elliptic curve
with Cremona label 15A8, with rank $0$ and precisely four
rational points. On the model $D_1$ the four
points are the two points at infinity,
and the points $(0,0)$ and $(1,0)$. 
None of these are images of real quadratic points on $X$.
We can conclude 
 that $u=\pm \sqrt{h_1(w)}$ is not rational. As
$h_1(w) \in \Q$ and $u \in K=\Q(\sqrt{d})$ we
have $h_1(w)=a^2 d$ for some non-zero $a \in \Q$. Thus
\[
h_2(w)=\frac{w^2 h_3(w)}{h_1(w)}=\frac{w^2 d}{a^2 d}=\frac{w^2}{a^2}.
\]
Therefore $v \in \Q$. 
Hence, if $P=(w,u,v)$ is a non-cuspidal
real quadratic point on $X$ then $\phi_2(P)=(w,v) \in D_2(\Q)$,
proving our claim.

To prove that the non-cuspidal real quadratic points 
on $X$ are modular we need an expression for the $j$-map.
The parametrization of $X$ via $g_1,\ldots,g_4$ allows us
to write down the $j$-map as a rational function on $X$
using the same method as described in the proof of Lemma~\ref{lem:b3s5}.
On the model \eqref{eqn:s3s5}, the $j$-map is given by
\begin{equation}\label{eqn:js3s5}
j=
\frac{(R_1+R_2\cdot v)^3}{8\, R_3^{15}}
\end{equation}
where
\begin{align*}
R_1 =&    w^{15} - 10 w^{14} + 45 w^{13} - 
	115 w^{12} + 155 w^{11} - 54 w^{10} + 185 w^9 -
        2395 w^8 + 10465 w^7\\
	& - 27055 w^6 + 47072 w^5 - 56090 w^4 + 43405 w^3 -
        19485 w^2 + 4020 w - 132;\\
R_2 =&    w^{13} - 8 w^{12} + 25 w^{11} - 35 w^{10} + 5 w^9 + 94 w^8 - 467 w^7 + 1805 w^6 -
        5175 w^5\\
	& + 9685 w^4 - 11216 w^3 + 7698 w^2 - 2715 w + 315;\\
R_3 =& w^2-w-1.
\end{align*}
The expression for $j$ involves only $w$ and $v$. Thus
the map $j: X \rightarrow X(1)$ factors via
$\phi_2: X \rightarrow D_2$. 
If $P$ is a non-cuspidal real quadratic point on $X$
then $\phi_2(P) \in D_2(\Q)$, and so $j(P) \in \Q$.
The modularity of $P$ follows from the modularity of elliptic curves over
the rationals \cite{modularity}.
\end{proof}

\noindent {\bf Remark}.
The reader might be surprised 
by the fact---crucial to our proof---that the $j$-map \eqref{eqn:js3s5}
on $X$ factors via the quotient $D_2$. 
Also surprising is the fact that the expression for $j$
is a perfect cube.
To explain this, we start
by the observation that
the normalizer $C_{\mathrm{s}}^+(3)$ 
of a split Cartan subgroup 
in $\GL_2(\F_3)$ is conjugate to a subgroup of the normalizer $C_{\mathrm{ns}}^+(3)$
of a non-split Cartan subgroup. Thus the map $X(\mathrm{s}3) \rightarrow X(1)$
factors via $X(\mathrm{ns}3)$. Moreover,  
$C_{\mathrm{s}}^+(3) \cap \SL_2(\F_3)$ has index $2$
in $C_{\mathrm{ns}}^+(3) \cap \SL_2(\F_3)$. It follows that the map
$X(\mathrm{s}3) \rightarrow X(\mathrm{ns}3)$ has degree $2$. Thus
the map $X \rightarrow X(1)$ factors via a degree $2$ map
$X \rightarrow X(\mathrm{ns}3,\mathrm{s}5)$. 
The curves $X(\mathrm{ns}3)$ and $X(\mathrm{s}5)$ have genus $0$.
In \cite[Section 5.1]{Chen}, Chen derives the following expressions
for the $j$-maps on $X(\mathrm{ns}3)$ and $X(\mathrm{s}5)$,
\[
j=n^3, \qquad 
j=\frac{(s^2-5)^3(s^2+5s+10)^3(s+5)^3}{(s^2+5s+5)^5} \, ,
\]
where $n$ and $s$ are respectively Hauptmoduln on $X(\mathrm{ns}3)$
and $X(\mathrm{s}5)$. Since $X \rightarrow X(1)$ factors
via $X(\mathrm{ns}3)$, it is clear that $j$ must be a cube on $X$,
and this provides a useful check on
our computation of $j$. A (singular) model for $X(\mathrm{ns}3,\mathrm{s}5)$
is given by
\[
n^3=\frac{(s^2-5)^3(s^2+5s+10)^3(s+5)^3}{(s^2+5s+5)^5} \, .
\] 
It is easy to see that this is birational to the elliptic curve
\[
E \; : \; s^2+5s+5= m^3.
\] 
The curve $E$ has Cremona label 225A1, and is in fact isomorphic to
$D_2$. This explains why $X \rightarrow X(1)$ factors through $D_2$.

\section{Models for some Modular Curves of Genus $1$}\label{sec:models}
We shall shortly study the quadratic points on
the four modular curves in \eqref{eqn:last4}.
\begin{lem}\label{lem:genera}
The four modular curves in \eqref{eqn:last4}  
respectively have genus $97$, $153$, $73$ and $113$.
\end{lem}
The lemma
follows from a standard application of Riemann--Hurwitz
applied to $X \rightarrow X(1)$, where $X$ runs through the 
four modular curves. For this, we need
the ramification indices of points above $0$, $1728$ 
and $\infty$, and we explain how to compute these in
Section~\ref{subsec:fb}. We omit the details. 

In view of the
very great genera, 
it is unlikely that methods 
used to compute the quadratic points on 
the curves in~\eqref{eqn:first3}
will succeed for the curves in~\eqref{eqn:last4}. 
Indeed, it is entirely impractical to write down smooth projective models
for any of these four curves. Instead, we shall
view them as normalizations of the fibre products
\begin{align*}
& X(\mathrm{b}3,\mathrm{b}5) \times_{X(1)} X(\mathrm{d}7),
& X(\mathrm{s}3,\mathrm{b}5) \times_{X(1)} X(\mathrm{d}7),\\
& X(\mathrm{b}3,\mathrm{b}5) \times_{X(1)} X(\mathrm{e}7),
& X(\mathrm{s}3,\mathrm{b}5) \times_{X(1)} X(\mathrm{e}7).
\end{align*}
We therefore write down models for 
\begin{equation}\label{eqn:four}
X(\mathrm{b}3,\mathrm{b}5), \qquad
X(\mathrm{s}3,\mathrm{b}5), \qquad
X(\mathrm{d}7), \qquad
X(\mathrm{e}7)
\end{equation}
 and the corresponding $j$-maps to $X(1)$.
We will also compute the Mordell--Weil groups 
of their Jacobians over $\Q$. We then
use this information to study the quadratic points
on the four curves in \eqref{eqn:last4}.
The  
curves in \eqref{eqn:four} are all elliptic curves (genus $1$ 
curves with $\Q$-points). Rather than identify elements of the
Mordell--Weil groups with points on the curves (as is usual
for elliptic curves), we shall view them as elements of $\Pic^0$;
we have found this to be more convenient for studying quadratic points. 

\subsection{The Modular Curves $X(\mathrm{d}7)$ and $X(\mathrm{e}7)$}
\begin{lem}\label{lem:de7}
A model for the curve $X(\mathrm{d}7)$ is given by 
\[
X(\mathrm{d}7)\; :\; y^2= -7 (x^4 - 10 x^3 + 27 x^2 - 10 x - 27).
\]
The $j$-map $X(\mathrm{d}7) \rightarrow X(1)$ is given by
\begin{equation}\label{eqn:xs7}
j=
\frac{
x
(x+1)^3 
(x^2-5 x+8)^3 (x^2-5 x+1)^3 (x^4-5 x^3+8 x^2-7 x+7)^3 
}
{(x^3-4 x^2+3 x+1)^7} \, .
\end{equation}
The $\Q$-points on $X(\mathrm{d}7)$ are $(5/2,\pm 7/4)$. The Mordell--Weil
group over $\Q$ is given by
\[
\frac{\Z}{2\Z} \cdot [ (5/2,-7/4) - (5/2,7/4)  ] \, .
\]
A model for the curve $X(\mathrm{e}7)$ is given by
\[
X(\mathrm{e}7) \; : \; y^2=7 (16 x^4 + 68 x^3 + 111 x^2 + 62 x + 11).
\]
The $j$-map $X(\mathrm{e}7) \rightarrow X(1)$ is given by
\begin{equation}\label{eqn:xns7}
j=\frac{
(3 x+1)^3
(4 x^2+5 x+2)^3 (x^2+3 x+4)^3 (x^2+10 x+4)^3 
}
{(x^3+x^2-2 x-1)^7} \, .
\end{equation}
The $\Q$-points on $X(\mathrm{e}7)$ are $(-1/3,\pm 14/9)$. The Mordell--Weil
group over $\Q$ is given by
\[
\frac{\Z}{2\Z} \cdot [ (-1/3,-14/9) - (-1/3,14/9)  ] \, .
\]
\end{lem}
\noindent \textbf{Remarks.}
\begin{itemize}
\item The models given for $X(\mathrm{d}7)$ and $X(\mathrm{e}7)$
are \lq quartic\rq\ genus $1$ curves. 
As both have rational points, we can also give
Weierstrass models. In fact, a Weierstrass model
for $X(\mathrm{d}7)$ is
\[
 y^2 + x y = x^3 - x^2 - 107 x + 552
\]  
with Cremona label 49A3. A Weierstrass model
for $X(\mathrm{e}7)$ is
\[
y^2 + x y = x^3 - x^2 - 1822 x + 30393
\]
with Cremona label 49A4. We prefer to work
with the models in the lemma, as the 
$j$-maps are expressible simply in terms of
the $x$-coordinates, a fact that is crucial
for the proof of Lemma~\ref{lem:b3b5d7}.
\item A model for $X(\mathrm{d}7)$ 
together with the $j$-map has been computed
by Elkies and is quoted 
 by Sutherland \cite[Section 3]{Sutherland}.
Elkies' model and $j$-map are equivalent to ours. 
\end{itemize}

\begin{proof}[Proof of Lemma~\ref{lem:de7}]
Recall that $C_\mathrm{s}^+(7)$ 
and $C_{\mathrm{ns}}^+(7)$
respectively denote the 
normalizers of the split and non-split Cartan subgroups in $\GL_2(\F_7)$.
By definition,
curves $X(\mathrm{d}7)$ and $X(\mathrm{e}7)$ are in fact $X(H_1)$ and $X(H_2)$
where $H_1$, $H_2$ are as in Proposition~\ref{prop:large2}.
Then $H_i \cap \SL_2(\F_7)$ with $i=1$, $2$ 
is a subgroup of index $2$
respectively in $C_\mathrm{s}^+(7) \cap \SL_2(\F_7)$ and $C_\mathrm{ns}^+(7) \cap
\SL_2(\F_7)$.  
Thus the inclusion $H_2 \subset C_{\mathrm{ns}}^+(7)$ induces
the following commutative diagram:
\[
\xymatrixcolsep{5pc}
\xymatrix{
X(\mathrm{e}7) \ar[r]^{\phi} \ar[rd]_j  & X(\mathrm{ns}7)  \ar[d]^j \\
& X(1)
}
\]
where $\phi$ has degree $2$. We also obtain a similar diagram
in which 
$X(\mathrm{d}7)$ 
and $X(\mathrm{s}7)$ respectively
replace $X(\mathrm{e}7)$ and $X(\mathrm{ns}7)$. 
The curves $X(\mathrm{s}7)$ and $X(\mathrm{ns}7)$
are rational.
In
\cite[Section 5.1]{Chen}, Chen derives expressions for
the maps $j: X(\mathrm{s}7) \rightarrow X(1)$
and $j : X(\mathrm{ns}7) \rightarrow X(1)$. These 
are given above in \eqref{eqn:xs7} and \eqref{eqn:xns7}
where $x$
is respectively a Hauptmodul for $X(\mathrm{s}7)$, $X(\mathrm{ns}7)$.
We used these to derive the models given for $X(\mathrm{d}7)$ and $X(\mathrm{e}7)$.
We only give the details for $X(\mathrm{e}7)$; the details for $X(\mathrm{d}7)$
are an easy adaptation of what follows.

From \eqref{eqn:xns7}, the map
$X(\mathrm{ns}7) \rightarrow X(1)$ has degree $21$. 
Subtracting $1728$ from both sides of \eqref{eqn:xns7}
and factoring gives
\begin{equation}\label{eqn:jm1728}
j-1728=
\frac{
    (16 x^4 + 68 x^3 + 111 x^2 + 62 x + 11)\cdot G(x)^2
}
{(x^3+x^2-2 x-1)^7},
\end{equation}
where
\[
G(x)=
7 (9 x^4 + 26 x^3 + 34 x^2 + 20 x + 4)
        (x^4 + 6 x^3 + 17 x^2 + 10 x + 2).
\]
It is convenient to view both $j: X(\mathrm{ns}7) \rightarrow X(1)$
and $j: X(\mathrm{e}7) \rightarrow X(1)$
as Belyi covers  
branched at $0$, $1728$ and $\infty$. 
The factorizations in \eqref{eqn:xns7}
and \eqref{eqn:jm1728} yield the ramification data
for $X(\mathrm{ns}7) \rightarrow X(1)$:
\begin{itemize} 
\item there are $7$ points above $0$, say $p_1, \dots, p_{7}$, 
all having ramification degree $3$;
\item there are $13$ points above $1728$, 
say $q_1,\dots, q_{13}$, with $q_1,\dots,q_5$ unramified and
$q_6,\dots,q_{13}$ having ramification degree $2$;
\item there are $3$ points above $\infty$, say
$r_1$, $r_2$, $r_3$, all having ramification degree $7$.
\end{itemize}
The corresponding ramification data for 
$X(\mathrm{e}7) \rightarrow X(1)$ 
can be deduced from the group $H_2$; for this
we applied the formulae
of Ligozat \cite{Ligozat}, which give the following information:
\begin{itemize}
\item there are $14$ points above $0$, say $P_1, \dots, P_{14}$, 
all having ramification degree $3$;
\item there are $22$ points above $1728$, 
say $Q_1,\dots, Q_{22}$,
with $Q_1$, $Q_2$ unramified, and $Q_{3},\dots, Q_{22}$
having ramification degree $2$;
\item there are $6$ points above $\infty$, say
$R_1,\dots,R_6$, all having ramification degree $7$.
\end{itemize}
It is clear from the above data that the 
degree $2$ cover $\phi : X(\mathrm{e}7) \rightarrow X(\mathrm{ns}7)$
is branched precisely at four of the five points
$q_1,\dots,q_5$ and nowhere else.
One of these five points is $\infty \in X(\mathrm{ns}7)(\Q)$
and the other four correspond to the roots of the irreducible factor
$16 x^4 + 68 x^3 + 111 x^2 + 62 x + 11$ in \eqref{eqn:jm1728}.
Since the set of branch points is stable under the action of 
Galois, we see that $\phi$ is branched precisely at the
roots of $16 x^4 + 68 x^3 + 111 x^2 + 62 x + 11$.
It follows that $X(\mathrm{e}7)$ has an equation of the form
\[
D_c~:~y^2= c (16 x^4 + 68 x^3 + 111 x^2 + 62 x + 11 ),
\]
where $c$ is a non-zero, squarefree integer. 
The map $\phi$ is then given by $(x,y) \mapsto x$,
and the map $j: X(\mathrm{e}7) \rightarrow X(1)$ is given
by the expression in \eqref{eqn:xns7}.
As $H_2$ is a congruence group of level $7$,
the modular curve $X(\mathrm{e}7)$ has a model
over $\Spec(\Z[1/7])$, and so its Jacobian
elliptic curve has good reduction outside $7$.
Using the formulae in \cite{Fisher},
the Jacobian elliptic curve of $D_c$ 
is given by
\[
E_c~:~Y^2 = X^3 + 111 c X^2 + 3512 c^2 x + 34224.
\]
The model $E_c$ has discriminant $2^{12} \cdot 7^3 \cdot c^6$.
If $p \mid c$ is a prime $\ne 7$ then the minimal model for 
$E_c$ has bad reduction at $p$.
Thus $c$ is not divisible by any prime $p \ne 7$.
It follows that $c=\pm 1$, $\pm 7$.
We checked that the only values for which the minimal model of $E_c$ 
has good reduction at $2$ are $c=-1$ and $c=7$.
Thus one of the two curves $D_{-1}$,  $D_{7}$ is a model for $X(\mathrm{e}7)$.
The curve $D_{-1}$ has no real points. It is however easy to see
that $X(\mathrm{e}7)$ must have real points. For example
if $E/\R$ is any real elliptic curve, then 
the image of $\overline{\rho}_{E/\R,7}$ has order $2$ and every
subgroup of $\GL_2(\F_7)$ of order $2$ is conjugate to a subgroup
of $H_2$, thus $E/\R$ gives rise to a real point on $X(\mathrm{e}7)$. 
It follows that $D_7$ is a model for $X(\mathrm{e}7)$.

This model has two obvious rational points $(-1/3,\pm 14/9)$
(both with $j$-invariant $0$). Taking $X(\mathrm{e}7)$ into Weierstrass form
gives the elliptic curve with Cremona label 49A4,
which has Mordell--Weil group isomorphic to $\Z/2\Z$.
Thus $(-1/3, \pm 14/9)$
are all the rational points of $X(\mathrm{e}7)$.
Finally the 
divisor $(-1/3,-14/9)-(-1/3,14/9)$ cannot be principal
as $X(\mathrm{e}7)$ has genus $1$, and so no two points
are linearly equivalent. Hence it must
generate $\Pic^0 X(\mathrm{e}7)$.
\end{proof}

\subsection{The Modular Curves $X(\mathrm{b}3,\mathrm{b}5)$ and $X(\mathrm{s}3,\mathrm{b}5)$}

In Lemmas~\ref{lem:b3b5} and~\ref{lem:s3b5} below we give
models for $X(\mathrm{b}3,\mathrm{b}5)$ and $X(\mathrm{s}3,\mathrm{b}5)$
as well as their rational points and the Mordell--Weil groups of their
Jacobians. 
In fact most of this is known (for example
\cite[Section 2]{Rubin}) but we know of no reference for the
$j$-maps, which are essential for our later proofs.
The method of Galbraith (sketched in the proof of 
Lemma~\ref{lem:b3s5}) applies only to non-hyperelliptic
modular curves of genus $\ge 3$, for which weight $2$
cusp forms can be used to define an embedding into projective
space. We now sketch an alternative
based on a paper of Mui\'{c} \cite[Theorem 3-3]{Muic}. 
Let $X$ be a modular curve
of genus $g$ corresponding to a congruence subgroup 
$\Gamma \subseteq \SL_2(\Z)$. Let $k \ge 3$ and 
let $S_k(\Gamma)$ be the space of cusp forms for $\Gamma$
of weight $k$. Write $t$ for the dimension $S_k(\Gamma)$,
and let $f_1,\dots,f_t$ be a basis. If $t \ge g+2$ then
$f_1,\ldots,f_t$ define an embedding $\lvert \mathfrak{c} \rvert$ of $X$ into
$\PP^{t-1}$ given by a very ample divisor $\mathfrak{c}$ 
of degree $t+g-1$. In particular, if $t \ge g+3$,
and so $\deg(\mathfrak{c}) \ge 2g+2$, the image of the embedding
is a scheme-theoretic intersection of quadrics \cite[page 130]{DaS}. By using
higher weight cusp forms where the space has sufficiently
large dimension one can straightforwardly adapt the method sketched in
the proof of Lemma~\ref{lem:b3s5} to write down a projective
model for $X$ and the associated $j$-map $X \rightarrow X(1)$.

For $X(\mathrm{b}3,\mathrm{b}5)$ and $X(\mathrm{s}3,\mathrm{b}5)$
the spaces of cusp forms of weight $4$ have dimensions $4$ and $8$
respectively, and this allowed us to write models for these
curves as intersections of quadrics in $\PP^3$ and $\PP^7$
respectively, as well as expressions for the $j$-maps.
Transforming these into Weierstrass models 
gives Lemmas~\ref{lem:b3b5} and~\ref{lem:s3b5} below;
the expressions for the $j$-maps on these two curves
are complicated but can be found in the accompanying {\tt Magma}
script (see page \pageref{magma}).

\begin{lem}\label{lem:b3b5}
The curve $X(\mathrm{b}3,\mathrm{b}5)=X_0(15)$ is 
an elliptic curve with the Weierstrass model
\[
X(\mathrm{b}3,\mathrm{b}5) \; : \;
y^2 + x y + y = x^3 + x^2 - 10 x - 10
\]
and Cremona label 15A1. Its $\Q$-points are 
\[
 \infty,\;
 (-1 , 0 ),\;
 (-2 , -2 ),\;
 (8 , -27 ),\;
 (3 , -2 ),\;
 (-13/4 , 9/8 ),\;
 (-2 , 3 ),\;
 (8 , 18 ).
\]
The Mordell--Weil group over $\Q$ is 
\[
\frac{\Z}{4\Z} \cdot [(-1,0)-\infty] \oplus \frac{\Z}{2\Z} [(-2,-2)-\infty].
\]
\end{lem}

\begin{lem}\label{lem:s3b5}
The curve $X(\mathrm{s}3,\mathrm{b}5)$
is an elliptic curve with the Weierstrass model
\begin{equation}\label{eqn:s3b5}
X(\mathrm{s}3,\mathrm{b}5) \; : \;
y^2+xy+y=x^3+x^2-5x+2,
\end{equation}
and Cremona label 15A3. The $\Q$-points 
of this model are
\[
 \infty,\; (3/4 , -7/8 ),\; (0 , -2 ),\; (2 , -4 ),\; (1 , -1 ),\;
 (-3 , 1 ),\; (0 , 1 ),\; (2 , 1 ).
\] 
The Mordell--Weil group of the Jacobian over $\Q$ is
\[
\frac{\Z}{4\Z} \cdot 
[(0,-2)-\infty]
\oplus \frac{\Z}{2\Z} \cdot 
[(3/4,-7/8) - \infty].
\]
\end{lem}

\section{A Sieve for the Symmetric Square of a Fibre Product}\label{sec:sieve}
Let $X_1$, $X_2$ be curves of genus $\ge 1$ defined 
over $\Q$, and let $J_1$, $J_2$ be their Jacobians. 
Let $\mu_k \; : \; X_k \rightarrow \PP^1$
be morphisms defined over $\Q$ for $k=1$, $2$. 
Suppose that the fibre product
$X_1 \times_{\PP^1} X_2$ is irreducible
and let $X$ be its normalization.
Let $\pi_k : X \rightarrow X_k$
be the induced projection maps. 
Fix rational degree $2$ divisors $D_k$ on $X_k$ and let
\begin{equation}\label{eqn:AJ}
\iota_k \; : \; X^{(2)}_k \rightarrow J_k,
\qquad D \mapsto [D-D_k]
\end{equation}
be the corresponding Abel--Jacobi maps for the symmetric squares.
Let
\[
\alpha \; : \; X^{(2)} \rightarrow J_1 \times J_2,
\qquad
D \mapsto
\left(
\iota_1 (\pi_1 D),\,
\iota_2 ( \pi_2 D)
\right).
\]
Let $p_1,\ldots,p_r$ be primes of good reduction for $X$.
For each $p=p_i$ we have the commutative diagram
\[
\xymatrix{
X^{(2)}(\Q) \ar^{\alpha}[r] \ar_{\red}[d] & J_1(\Q) \times J_2 (\Q) \ar_{\red}[d]\\
X^{(2)}(\F_p) \ar^{\alpha}[r]  & J_1(\F_p) \times J_2 (\F_p) \\
}
\]
where the vertical maps denote reduction modulo $p$.
The following lemma is an obvious consequence of the diagram's commutativity.
\begin{lem}\label{lem:sieve}
With the above notation and assumptions,
\[
\alpha \left( X^{(2)}(\Q) \right)
\subseteq
\bigcap_{p=p_i} \red^{-1} \alpha \left( X^{(2)}(\F_p) \right).
\]
In particular, if the intersection on the right-hand side is empty,
then $X^{(2)}(\Q)=\emptyset$ and so the curve $X$ does not have any
rational or quadratic points.
\end{lem}

\bigskip
\noindent \textbf{Remark.} 
Let $X/\Q$ be a curve of genus $\ge 2$ 
and let $J$ be its Jacobian. Fix some $P_0 \in X(\Q)$
and let $\iota : X \rightarrow J$ be the corresponding
Abel--Jacobi map $P \mapsto [P-P_0]$. 
For any prime $p$ of good reduction for $X$, we have a commutative diagram
\[
\xymatrix{
X(\Q) \ar^{\iota}[r] \ar_{\red}[d] & J(\Q) \ar^{\red}[d] \\
X(\F_p) \ar^{\iota}[r]  & J(\F_p)  \\
}
\]
If $p_1,\dotsc,p_r$ we have 
\[
\iota \left( X(\Q) \right)
\subseteq
\bigcap_{p=p_i} \red^{-1} \iota \left( X(\F_p) \right).
\]
This sieve is known as the Mordell--Weil sieve (e.g. \cite[Section 3]{StollRat})
and gives information on the image of $X(\Q)$ inside $J(\Q)$. The Mordell--Weil
sieve can be extended in the obvious manner to the symmetric square of $X$
(see for example \cite[Section 5]{chabsym}), and then it gives
information on the image of $X^{(2)}(\Q)$ inside $J(\Q)$.
Now suppose $X$ is the normalization of a fibre product
$X_1 \times_{\PP^1} X_2$ as above. The Jacobian $J$ 
is in general much larger (in dimension) than $J_1 \times J_2$,
and the information provided by Lemma~\ref{lem:sieve}
will be weaker than that given by the Mordell--Weil sieve.
However, in the Section~\ref{sec:b3b5d7} where we apply Lemma~\ref{lem:sieve},
it is realistic to carry out explicit computations with
$J_1$, $J_2$ but not with $J$. Indeed, in our 
applications $J_1$ and $J_2$
will be elliptic curves, whereas the possible dimensions
for $J$ will be $97$, $153$, $73$ and $113$.

\subsection{Points on the Normalization of a Fibre Product} \label{subsec:fb}
To apply Lemma~\ref{lem:sieve}, we need to enumerate the points
of $X^{(2)}(\F_p)$ for several primes $p$. 
To enumerate $X^{(2)}(\F_p)$ it is enough to know the points
of $X(\F_p)$ and $X(\F_{p^2})$, and this is easy in practice if $p$
is sufficiently small and one has a convenient model
(say smooth and projective) for $X$.
For the four curves \eqref{eqn:last4} we will apply the lemma to,
it seems impractical to write
down such a model in view of the very large
genera of these curves (Lemma~\ref{lem:genera}). 
In what follows we will use a description of the points
of the normalization $X$ of $X_1 \times_{\PP^1} X_2$
in terms of the points of $X_1$, $X_2$. This will allow
us to enumerate the points of $X$ over finite fields.

We continue with the above notation, and we let $\F$
be either an extension of $\Q$ or an extension of $\F_p$ where $p$ is a prime of good reduction
for $X$, $X_1$, $X_2$ and the projection maps, and such that the
characteristic of $\F$ does not divide 
the ramification indices of the maps $\mu_k$.
By definition of $X_1 \times_{\PP^1} X_2$, 
if $P \in X(\F)$, and $P_k=\pi_k(P)$ for $k=1$, $2$
then $P_k \in X_k(\F)$ and $\mu_1(P_1)=\mu_2(P_2)$.
Now suppose $P_1 \in X_1(\F)$, $P_2 \in X_2(\F)$
satisfy $\mu_1(P_1)=\mu_2(P_2)$. We would like to 
describe the points $P \in X(\F)$ such that
$P_k=\pi_k(P)$ for $k=1$, $2$. Here we follow the
description of Flynn and Testa \cite[Section 2]{FlynnTesta}.
Write $\gamma=\mu_1(P_1)=\mu_2(P_2) \in \PP^1(\F)$.
Suppose first that $\gamma \ne \infty$, 
so we may think of $\gamma \in \mathbb{A}^1(\F)$.
We consider $\mu_1$, $\mu_2$ as elements of the completions
of the local rings of $X_k$ at $P_k$, which we denote by
$\widehat{\OO_k}$. 
Then
\[
\mu_1=\gamma+\epsilon s^m, \qquad 
\mu_2=\gamma+\delta t^n,
\]
where $s \in \widehat{\OO_1}$, $t \in \widehat{\OO_2}$
are elements of valuation $1$, the coefficients
$\epsilon$ and $\delta$ belong to $\F^*$, and
$m$, $n$ are respectively the ramification indices for 
$\mu_1$, $\mu_2$ at $P_1$, $P_2$.
A local 
equation for $X_1 \times_{\PP^1} X_2$ in the neighbourhood
of $(P_1,P_2)$ is given by
\[
\epsilon s^m
=\delta t^n .
\]
Of course, on this local equation, the point $(P_1,P_2)$
is represented by $(0,0)$. Now let $d=\gcd(m,n)$.
Over $\overline{\F}$, the 
local equation 
has $d$ components:
$\lambda s^{m^\prime}=t^{n^\prime}$
where $m^\prime=m/d$, $n^\prime=n/d$ and $\lambda$ runs through the
$d$-th roots of $\epsilon/\delta$. We normalize these components
by letting $N=m^\prime \cdot n^\prime$ and making the substitutions
\[
s=\lambda^a  r^{n^\prime},\qquad t=\lambda^b r^{m^\prime},
\]
with $a$, $b \in \Z$ chosen so that $a m^\prime-b n^\prime=-1$.
It follows that there are precisely $d$ points $P \in X$
above $(P_1,P_2)$. We may represent these by $(P_1,P_2,\lambda)$
where $\lambda^d=\epsilon/\delta$. Such a point belongs
to $X(\F)$ if and only if $P_1 \in X_1(\F)$, $P_2 \in X_2(\F)$
and $\lambda \in \F$. Moreover, the parameter $r$ has valuation
$1$ in the completion of the local ring of $X$ at $P$.
The ramification index for $X \rightarrow \PP^1$
at each of these $P$ is 
$m n^\prime=m^\prime n= m n / \gcd(m,n)$.
To see this, note that 
the composition of maps
\[
X \rightarrow X_1 \times_{\PP^1} X_2 \rightarrow X_1 \rightarrow \PP^1
\]
viewed
as an element of the completion of the local ring of $X$ at $P$
can be written as 
\[
\mu_1= \gamma+ \epsilon s^m=
\gamma + \epsilon (\lambda^a r^{n^\prime})^m
=\gamma+\epsilon \cdot \lambda^a \cdot r^{m n^\prime}.
\]
The above argument needs an easy adaptation if $\gamma=\infty \in \PP^{1}$.
When representing the points of $X$ as $P=(P_1,P_2,\lambda)$
the projection maps are obvious: $\pi_k(P)=P_k$.

If $\F=\F_p$ or $\F_{p^2}$, the above gives an algorithm for constructing the
$\F$-points of $X$ from the $\F$-points of $X_1$, $X_2$ and so it enables
us to write down $X^{(2)}(\F_p)$.

\section{Real Quadratic Points on 
$X(\mathrm{b}3,\mathrm{b}5,\mathrm{d}7)$,
$X(\mathrm{s}3,\mathrm{b}5,\mathrm{d}7)$,
$X(\mathrm{b}3,\mathrm{b}5,\mathrm{e}7)$, 
and
$X(\mathrm{s}3,\mathrm{b}5,\mathrm{e}7)$.
}\label{sec:b3b5d7}

\begin{lem}\label{lem:b3b5d7}
The curve $X(\mathrm{b}3,\mathrm{b}5,\mathrm{e}7)$
does not have
any rational or quadratic points.
The non-cuspidal real quadratic points on the curves
\begin{equation}\label{eqn:3mod}
X(\mathrm{b}3,\mathrm{b}5,\mathrm{d}7),\qquad
X(\mathrm{s}3,\mathrm{b}5,\mathrm{d}7), \qquad
X(\mathrm{s}3,\mathrm{b}5,\mathrm{e}7) 
\end{equation}
are modular.
\end{lem}
\begin{proof}
Let $X=X(\mathrm{b}3,\mathrm{b}5,\mathrm{d}7)$ which
we think of as the normalization of $X_1 \times_{\PP^1} X_2$
where $X_1=X(\mathrm{b}3,\mathrm{b}5)$, $X_2=X(\mathrm{d}7)$
with the maps $\mu_k : X_k \rightarrow \PP^1$
being just the $j$-maps to $X(1)$. Write $J_1$, $J_2$
for the Jacobians of $X_1$, $X_2$.
In Section~\ref{sec:models} we determined models for the
genus $1$ curves $X_1$, $X_2$ as well as the Mordell--Weil
groups $J_1(\Q)$, $J_2(\Q)$. 
We fix degree $2$ divisors
\begin{equation}\label{eqn:D}
D_1=2 \infty, \qquad
D_2=2 \cdot (5/2,7/4)
\end{equation}
on $X_1$, $X_2$ respectively, which we use
to define Abel--Jacobi maps $X_k^{(2)} \rightarrow J_k$
as in \eqref{eqn:AJ}. We now apply Lemma~\ref{lem:sieve}
with $11 \leq p < 100$. We note that these primes are 
of good reduction for $X$, and that they do not divide the
ramification indices for the maps $X_k \rightarrow X(1)$.
We find that
\begin{equation}\label{eqn:sieveb3b5}
\begin{split}
\alpha \left( X^{(2)} (\Q) \right)
& \subseteq
\bigg\{ \; \Big( [0], \; [(5/2,-7/4)-(5/2,7/4)] \Big),  \\ 
& \qquad \Big([(3,-2)-\infty  ],\; [(5/2,-7/4)-(5/2,7/4)] \Big), \\
& \qquad
 \Big([(8,18)-\infty], \; [(5/2,-7/4)-(5/2,7/4)]\Big) \; \bigg\}\\ 
& \subseteq J_1(\Q) \times J_2(\Q).
\end{split}
\end{equation}
Now let $P$ be a quadratic point on $X$ and let $P^\sigma$ be its 
conjugate. Then $D=P+P^\sigma \in X^{(2)}(\Q)$. 
From \eqref{eqn:sieveb3b5}
and the definition of $\alpha$,
we have 
\[
\pi_2(P)+\pi_2(P^\sigma)-D_2 \; \thicksim \; (5/2,-7/4)-(5/2,7/4).
\]
Thus 
\[
\pi_2(P)+\pi_2(P^\sigma)=\divisor(f)+(5/2,7/4)+(5/2,-7/4)
\]
for some non-zero function $f$ on $X_2=X(\mathrm{d}7)$. The function $f$ belongs to
the Riemann--Roch space of $(5/2,7/4)+(5/2,-7/4)$.
A $\Q$-basis for this Riemann--Roch space is 
\[
\frac{1}{2x-5}, \qquad \frac{x}{2x-5}. 
\]
Thus
\[
f=\frac{a+b x}{2x-5}
\]
for some $a$, $b \in \Q$, not both zero.
Thus the
$x$-coordinate of $\pi_2(P)$ and $\pi_2(P^\sigma)$ is 
$-a/b \in \PP^{1}(\Q)$.
From Lemma~\ref{lem:de7}, we know that the $j$-value
of $\pi_2(P)$ is expressible in terms of its $x$-coordinate
as given by \eqref{eqn:xs7}, hence it belongs to $\PP^{1}(\Q)$.
Now suppose $P$ is a non-cuspidal real quadratic point on $X$.
It follows that any elliptic curve supporting $P$ has
$j$-invariant belonging to $\Q$ and so is modular by \cite{modularity}.

We apply the same sieve to $X=X(\mathrm{s}3,\mathrm{b}5,\mathrm{d}7)$,
where we take $X_1=X(\mathrm{s}3,\mathrm{b}5)$
and $X_2=X(\mathrm{d}7)$, and we take $D_1$, $D_2$
as in \eqref{eqn:D}, where of course, now $\infty$ is the point 
at infinity
on the model for $X_1$ given in \eqref{eqn:s3b5}.
We obtained
\begin{equation}\label{eqn:sieves3b5}
\begin{split}
\alpha \left( X^{(2)} (\Q) \right)
& \subseteq
\bigg\{ \; \Big( [0], \; [(5/2,-7/4)-(5/2,7/4)] \Big),  \\ 
& \qquad  \Big([(2,-4)-\infty  ],\; [(5/2,-7/4)-(5/2,7/4)] \Big), \\
& \qquad 
 \Big([(1,-1)-\infty], \; [(5/2,-7/4)-(5/2,7/4)]\Big) \; \bigg\}\\ 
& \subseteq J_1(\Q) \times J_2(\Q).
\end{split}
\end{equation}
By exactly the same argument as before, 
quadratic points have rational $j$-invariants,
and thus non-cuspidal real
quadratic points are modular.

Applying the sieve to $X=X(\mathrm{b}3,\mathrm{b}5,\mathrm{e}7)$
shows
that $\alpha(X^{(2)}(\Q))$ is empty, so these this curve have no
rational or quadratic points. 

Finally, we apply the sieve to 
$X=X(\mathrm{s}3,\mathrm{b}5,\mathrm{e}7)$.
Of course, we take $X_1=X(\mathrm{s}3,\mathrm{b}5)$, $X_2=X(\mathrm{e}7)$.
We fix the degree $2$ divisors
\[
D_1=2 \infty, \qquad
D_2=2 \cdot (-1/3,14/9)
\]
of $X_1$ and $X_2$. The sieve yields
\begin{equation*}
\begin{split}
\alpha \left( X^{(2)} (\Q) \right)
 \subseteq
\bigg\{ \; \Big( [(0,1)-\infty], \; [(-1/3,-14/9)-(-1/3,14/9)] \Big)
 \; \bigg\} .
\end{split}
\end{equation*}
If $P$ is a quadratic point on $X$ and $\pi_2(P)$ is its
image on $X_2=X(\mathrm{e}7)$ then
\[
\pi_2(P)+\pi_2(P^\sigma)=\divisor(f)+
(-1/3,14/9)+(-1/3,-14/9).
\]
As before,
the $x$-coordinate of $\pi_2(P)$ on the model for $X(\mathrm{e}7)$ 
belongs to $\PP^1(\Q)$ and so the $j$-invariant (given
by \eqref{eqn:xns7}) is rational. This completes the proof.
\end{proof}

\bigskip

\noindent \textbf{Remarks}.

\medskip

\noindent \textbf{(I)}
In the above, we have morphisms $\alpha : X^{(2)} \rightarrow J_1 \times J_2$.
The points $X^{(2)}(\Q)$ belong to the fibres $\alpha^{-1} (J_1(\Q) \times J_2(\Q))$.
Since $J_1(\Q)$ and $J_2(\Q)$ are finite and explicitly known, the reader
is perhaps puzzled as to why we do not use this to recover $X^{(2)}(\Q)$. 
We are in fact able to write down equations for
these fibres as schemes in $\PP^1 \times \PP^1$. Most of these
are zero-dimensional, but are of such high degree that 
we have found it impractical to decompose them using Gr\"{o}bner
basis computations. A few of the fibres are however $1$-dimensional.
These are so complicated that we are unable
to say anything about their rational points. 

\medskip

\noindent \textbf{(II)} Our sieve cannot 
prove that $X^{(2)}(\Q)$ is empty for 
the three modular curves \eqref{eqn:3mod}.
It is natural to ask whether these curves have rational or
quadratic points. The
genera of these three modular curves are $97$, $153$ 
and $113$ respectively
by Lemma~\ref{lem:genera}. As these are so large, it is reasonable 
to suspect that the only quadratic points lie above $0$, $1728$ or
$\infty$ (and Remark (III) below shows that this is indeed the case). 
Using the description of the points of 
$X=X(\mathrm{b}3,\mathrm{b}5,\mathrm{d}7)$
as the fibre product of $X_1=X(\mathrm{b}3,\mathrm{b}5)$
and $X_2=X(\mathrm{d}7)$ we find that there are precisely
$8$ quadratic cusps, and no quadratic points above $0$ and $1728$.
We shall denote the quadratic cusps by $P_1^{\pm}$, $P_2^{\pm}$,
$P_3^{\pm}$, $P_4^{\pm}$. These respectively lie above
\[
\big(\infty, \infty_{\pm}\big),\;
\big( (-1,0), \infty_{\pm} \big),\;
\big( (-2,3), \infty_{\pm} \big),\;
\big( (8,18), \infty_{\pm} \big) \;
 \in X_1 \times_{\PP^1} X_2.
\]
Here $\infty_{+}$, $\infty_{-} \in X_2$ are
conjugates defined over $\Q(\sqrt{-7})$.
Thus $P_i^{+}+P_i^{-}$ belong to $X^{(2)}(\Q)$. We find that
\begin{align*}
& \alpha(P_1^{+}+P_1^{-})=
\alpha(P_2^{+}+P_2^{-})=
\big( \, [0],\; [ (5/2,-7/4)- (5/2,7/4)] \, \big),\\
& \alpha(P_3^{+}+P_3^{-})=
\alpha(P_4^{+}+P_4^{-})=
\big(\, [(3,-2)-\infty],\; [(5/2,-7/4)-(5/2,7/4)] \, \big).
\end{align*}
This is consistent with the result of the sieve given in \eqref{eqn:sieveb3b5},
and provides a useful check on the correctness of our computations.

We carried the same computation for $X=X(\mathrm{s}3,\mathrm{b}5,\mathrm{d}7)$
and found quadratic cusps $P_1^{\pm}$, $P_2^{\pm}$,
$P_3^{\pm}$, $P_4^{\pm}$ and no other quadratic points
above $0$, $1728$ and $\infty$. The four quadratic cusps respectively lie above
\[
\big(\infty, \infty_{\pm}\big),\;
\big( (0,1), \infty_{\pm} \big),\;
\big( (2,-4), \infty_{\pm} \big),\;
\big( (-3,1), \infty_{\pm} \big) \;
 \in X_1 \times_{\PP^1} X_2,
\]
where of course $X_1=X(\mathrm{s}3,\mathrm{b}5)$ and $X_2=X(\mathrm{d}7)$.
We find that
\begin{align*}
&\alpha(P_1^{+}+P_1^{-})=
\alpha(P_4^{+}+P_4^{-})=
\big(\, [0],\; [(5/2,-7/4)-(5/2,7/4)]\, \big),\\
&\alpha(P_2^{+}+P_2^{-})=
\alpha(P_3^{+}+P_3^{-})=
\big(\, [(1,-1)-\infty],\; [(5/2,-7/4)-(5/2,7/4)]\, \big).
\end{align*}
This is  
consistent with the result of the sieve given in \eqref{eqn:sieves3b5}.

As a final check on our computations, we verified
that there are no quadratic points 
on $X(\mathrm{b}3,\mathrm{b}5,\mathrm{e}7)$
and $X(\mathrm{s}3,\mathrm{b}5,\mathrm{e}7)$ 
above $0$, $1728$ and $\infty$.

\medskip

\noindent \textbf{(III)}
Although we do not need this for the proof of Lemma~\ref{lem:b3b5d7}, it is interesting to
ask if it can be shown that all quadratic points $P$ 
on $X_1:=X(\mathrm{b}3,\mathrm{b}5,\mathrm{d}7)$,
$X_2:=X(\mathrm{s}3,\mathrm{b}5,\mathrm{d}7)$,
$X_3:=X(\mathrm{s}3,\mathrm{b}5,\mathrm{e}7)$ do satisfy $j(P) \in \{0,1728,\infty\}$.
We explain how to solve this problem  by reducing it to the determination of
{\em rational} points on some modular curves with \lq milder\rq\ level structure.
Indeed, let $P$ be a non-cuspidal quadratic point on $X$, where $X$ is one
of these three modular curves $X_1$, $X_2$, $X_3$. We know by the proof of Lemma~\ref{lem:b3b5d7}
that $j(P) \in \Q$. Thus the point $P$ is supported by some elliptic curve $E$
defined over $\Q$, even though the level structure is defined over some quadratic
field $K$. We know in particular that $\overline{\rho}_{E,5}(G_K)$ is contained
in a Borel subgroup $B(5)$ of $\GL_2(\F_5)$. Now $\overline{\rho}_{E,5}(G_K)$
is a subgroup of $\overline{\rho}_{E,5}(G_\Q)$ of index $1$ or $2$. By a straightforward
analysis of the subgroups of $\GL_2(\F_5)$ it follows that $\overline{\rho}_{E,5}(G_\Q)$
is contained either in $B(5)$ or $C_{\mathrm{s}}^{+}(5)$ (the normalizer of a
split Cartan subgroup). Similarly, for $X=X_1$, we have $\overline{\rho}_{E,3}(G_\Q)$
is contained either in $B(3)$ or $C_{\mathrm{s}}^{+}(3)$ and for $X=X_2$, $X=X_3$
it is contained in $C_{\mathrm{ns}}^{+}(3)$.  Finally, for $X=X_1$, $X_2$,
the image $\overline{\rho}_{E,7}(G_\Q)$ is contained in $C_{\mathrm{s}}^{+}(7)$ and
for $X=X_3$ it is contained in $C_{\mathrm{ns}}^{+}(7)$. It is now not 
difficult to show in all cases, with the help of the results of this paper and of \cite{Bao}, that $j(E) \in \{0,1728\}$.
Consider $X=X_1$ for example. Then 
$E$ gives rise to a $\Q$-point of one of the curves
$X(\mathrm{b}3,\mathrm{b}5)$, $X(\mathrm{b}3,\mathrm{s}5)$,
$X(\mathrm{s}3,\mathrm{b}5)$ or $X(\mathrm{s}3,\mathrm{s}5)$. We have
already determined the $\Q$-points on these four curves in the proofs of 
Lemmas~\ref{lem:b3s5},~\ref{lem:s3s5}, and in Lemmas~\ref{lem:b3b5},~\ref{lem:s3b5}. 
Enumerating the $j$-invariants of these points $Q$, we are able to show for $j(Q) \ne 0$, $1728$
that they do not pull back to non-cuspidal quadratic points on $X(\mathrm{d}7)$.

\section{Modularity of Elliptic Curves over $\Q(\sqrt{5})$}\label{sec:Qsqrt5}
To complete the proof of Theorem~\ref{thm:modquadratic} we merely
have to prove Lemma~\ref{lem:Qsqrt5}, which states
that every elliptic curve over $K=\Q(\sqrt{5})$ is modular.
Again we reduce this problem to showing that the $K$-points
on certain modular curves are in fact modular (in
the sense that they correspond to modular elliptic curves).
We give fewer details in this section as compared
to previous sections. This is because the task
is much simpler: determining
the points on a curve over one quadratic field is usually 
easier than determining the points over all real quadratic fields.

\begin{proof}[Proof of Lemma~\ref{lem:Qsqrt5}]
Let $K=\Q(\sqrt{5})$, and let $E$ be an elliptic curve
over $K$. We would like to show that $E$ is modular. By
Theorems~\ref{thm:35} and~\ref{thm:7},
if $\overline{\rho}_{E,p}(G_{K(\zeta_p)})$ is absolutely
irreducible with $p=3$ or $7$ then $E$ is modular.
Thus we may restrict
our attention to elliptic curves $E$ such that
$\overline{\rho}_{E,p}(G_{K(\zeta_p)})$ is 
simultaneously absolutely reducible
for $p=3$ and $7$. 
Now appealing to Corollary~\ref{cor:large2} we see that
$E$ gives rise to a $K$-point on one of the following 
four modular curves
\[
X(\mathrm{d}7), \qquad X(\mathrm{e}7), \qquad X(\mathrm{b}3,\mathrm{b}7),
\qquad X(\mathrm{s}3,\mathrm{b}7).
\]
We determined models for the first two curves in Lemma~\ref{lem:de7}
and found that they are isomorphic (over $\Q$) to elliptic curves 49A3 and 49A4
respectively. The third curve is isomorphic (over $\Q$) to the elliptic curve
21A1. Let $X$ be any of the first three curves.
By considering
both $X$ and its $\sqrt{5}$-twist it is easy to show that $X(K)=X(\Q)$
in all cases. 
It follows that the $K$-points on the
first three curves are modular.

Finally let $X=X(\mathrm{s}3,\mathrm{b}7) \cong X_0(63)/\langle w_9 \rangle$.
This curve has the obvious quotient $Y=X_0(63)/\langle w_7, w_9 \rangle$.
The curve $Y$ has genus $1$ and is in fact isomorphic to the
elliptic curve 21A1. By the above, $Y(K)=Y(\Q)$.
Thus if $P$ is a non-cuspidal $K$-point on $X$ then its image
on $Y$ belongs to $Y(\Q)$. Let $\sigma$ be the non-trivial
automorphism of $K$. Thus $P$ and $P^\sigma$ map to the same
point on $Y$. Hence $P^\sigma=P$ or $P^\sigma=w_7(P)$. In the former
case $P$ is a rational point and hence modular. In the latter case,
$P$ corresponds to a $\Q$-curve, and is therefore modular. 
\end{proof} 

\part{Further Directions}
\section{On Modularity over Totally Real Fields (Again)}\label{sec:totreal}
Most of the theory we have developed in Parts 2, 3 and 4 applies to elliptic
curves over arbitrary totally real fields. If one is interested in
modularity of elliptic curves over a particular totally real
field $K$ then it might be possible to prove this by the methods of this paper, together with
explicit methods for computing rational points on curves (e.g. 
 \cite{chabnf},
\cite{StollRat}). 
However, one can considerably simplify this task by imposing
suitable restrictions on the field $K$ and
elliptic curve $E$.
\begin{thm}\label{thm:unram}
Let $p=5$ or $7$. Let $K$ be a totally real field
having some unramified place $\upsilon$ above $p$.
Let $E$ be an elliptic curve semistable at $\upsilon$ and suppose $\overline{\rho}_{E,p}$ 
is irreducible. Then
$E$ is modular.
\end{thm}
\begin{proof}
As $\upsilon/p$ is unramified, we have $K \cap \Q(\zeta_p)=\Q$.
Let $I_\upsilon$ be the inertia subgroup of $G_K$ at $\upsilon$.
By \cite[Propositions 11, 12, 13]{Serre}, the image of $\overline{\rho}_{E,p}(I_\upsilon)$
in $\PGL_2(\F_p)$ contains a cyclic subgroup of order $p-1$ or $p+1$. 
It follows from Proposition~\ref{prop:large2} that $\overline{\rho}(G_{K(\zeta_p)})$
is absolutely irreducible. Thus $E$ is modular by Theorems~\ref{thm:35}
and~\ref{thm:7}.
\end{proof}
In \cite{recipes}, the first-named author proves a version of Theorem~\ref{thm:unram} with $p=3$
(and without the irreducibility assumption) using modularity lifting theorems due to
Skinner and Wiles \cite{SWI}, \cite{SWII}.

\section{Diophantine Applications}\label{sec:dio}
Wiles' proof of Fermat's Last Theorem provides a template for 
attacking many other Diophantine equations via Frey curves and modularity (see \cite{BCDY},
\cite{Siksek} 
for recent surveys). 
Modularity results over totally real fields extend the applicability
of this template, as for example in \cite{DF2}. 
We expect the results and methods
of the current paper to play an important r\^ole in the development
of this area. Indeed, a forthcoming paper \cite{FS} uses Theorem~\ref{thm:modquadratic}
 to study the Fermat equation over real quadratic fields. In particular,
it is shown that if $d>6$ is a squarefree integer satisfying $d \equiv 3 \pmod{8}$
or $d \equiv 6$, $10 \pmod{16}$, then there is an effectively computable
constant $B_d$ such that for all primes $p > B_d$, every
 solution to the Fermat equation 
\[
x^p+y^p=z^p, \qquad x,y,z \in \Q(\sqrt{d})
\]
satisfies $xyz=0$.


\begin{thebibliography}{}


\bibitem{Allen}
P.\ B.\ Allen, {\em Modularity of nearly ordinary 2-adic residually
  dihedral Galois representations}, {\tt arXiv:1301.1113v2}, 11 September 2013.

\bibitem{BC}
B.\ S.\ Banwait and J.\ Cremona,
{\em Tetrahedral elliptic curves and the local-to-global 
principle for isogenies}, 
{\tt arXiv:1306.6818v2}, 16 August 2013. 


\bibitem{BGG1}
T.\ Barnet-Lamb, T.\ Gee and D.\ Geraghty,
{\em Congruences between Hilbert modular forms:
constructing ordinary lifts},
Duke Math. Journal {\bf 161} (2012), 1521--1580.

\bibitem{BGG2}
T.\ Barnet-Lamb, T.\ Gee and D.\ Geraghty,
{\em Congruences between Hilbert modular forms:
constructing ordinary lifts II},
Mathematical Research Letters {\bf 20} (2013), 81--86.


\bibitem{bars2} F.\ Bars, A.\ Kontogeorgis and X.\ Xarles,
{\em Bielliptic and hyperelliptic modular curves $X(n)$ and 
the group $\Aut(X(n))$}, Acta Arithmetica, to appear.






\bibitem{BCDY} M.\ A.\ Bennett, I.\ Chen, S.\ R.\ Dahmen and S.\ Yazdani, 
{\em Generalized Fermat equations: a miscellany}, preprint, 2013.


\bibitem{MAGMA} W.\ Bosma, J.\ Cannon and C.\ Playoust: {\em The Magma
Algebra System I: The User Language}, J.\ Symb.\ Comp.\ {\bf 24} (1997),
235--265. (See also {\tt http://magma.maths.usyd.edu.au/magma/})



\bibitem{Blasius-Rogawski}
D.\ Blasius and J.\ D.\ Rogawski, {\em Motives for Hilbert modular
  forms}, Inventiones Mathematicae {\bf 114} (1993), no.~1,
  55--87.

\bibitem{modularity} C.\ Breuil, B.\ Conrad, F.\ Diamond, R.\ Taylor,
 {\em On the modularity of elliptic curves over $\Q$: wild $3$-adic exercises},
 Journal of the American Mathematical Society {\bf 14} (2001), 843--939.

\bibitem{BreuilDiamond} C.\ Breuil and F.\ Diamond
{\em Formes modulaires de Hilbert modulo $p$ et valeurs d'extensions galoisiennes},
Annales Scientifiques de l'\'{E}cole Normale Sup\'{e}rieure, to appear.

\bibitem{BN} P.\ Bruin and F.\ Najman,
{\em On isogenies of elliptic curves over quadratic fields},
preprint, 2013.

\bibitem{Carayol1} H.\ Carayol, {\em Sur les repr\'{e}sentations $\ell$-adiques associ\'{e}es aux formes modulaires de Hilbert}, Ann.\ Sci.\ Ec.\ Norm.\ Sup.\ {\bf 4} (1986), 409--468.


\bibitem{Carayol2} H.\ Carayol, {\em Sur les repr\'{e}sentations galoisiennes modulo attach\'{e}es aux formes modulaires}, Duke Math.\ J.\ {\bf 59} (1989), 785--801.



\bibitem{Chen} I.\ Chen,
{\em The Jacobian of Modular Curves Associated 
to Cartan Subgroups},
DPhil thesis, University of Oxford, 1996.

\bibitem{CY} 
J.\ Coates and S.\ T.\ Yau (editors),
{\em Elliptic Curves, Modular Forms \& Fermat's Last Theorem},
second edition,
International Press,
  Cambridge, {MA}, 1997.

\bibitem{Cohen} H.\ Cohen,
{\em Number Theory, Volume II: Analytic and Modern Tools},
GTM {\bf 240}, Springer-Verlag, 2007.

\bibitem{CDT} B.\ Conrad,
F.\ Diamond and R.\ Taylor,
{\em Modularity of certain potentially Basotti--Tate
Galois representations},
Journal of the American Mathematical Society
{\bf 12} (1999), no.\ 2, 521--567.

\bibitem{CS} G.\ Cornell and J.\ H.\ Silverman (editors),
{\em Arithmetic Geometry},
Springer-Verlag, 1986.

\bibitem{Cornell} G.\ Cornell, J.\ H.\ Silverman
and G.\ Stevens (editors),
{\em Modular Forms and Fermat's Last Theorem},
Springer-Verlag, 1997.

\bibitem{Cre} J.\ E.\ Cremona,
{\em Algorithms for Modular Elliptic Curves},
Cambridge University Press, 1992.

\bibitem{modcurves} J.\ Cremona, J.-C.\ Lario, J.\ Quer, K.\ Ribet (editors),
{\em Modular Curves and Abelian Varieties}, Progr.\ Math.\ {\bf 224}, 
Birkh\"{a}user, Basel, 2004.

\bibitem{DaS} V.\ I.\ Danilov and V.\ V.\ Shokurov,
{\em Algebraic Curves, Algebraic Manifolds and Schemes},
Springer-Verlag, 1998.

\bibitem{DDT}
H.\ Darmon, F.\ Diamond, and R.\ Taylor, {\em Fermat's last theorem},
pages 2--140 of \cite{CY}.

\bibitem{Diamond} F.\ Diamond,
{\em On deformation rings and Hecke rings},
Annals of Mathematics {\bf 144} (1996), 137--166.


\bibitem{DS} F.\ Diamond and J.\ Shurman,
{\em A First Course on Modular Forms},
GTM {\bf 228}, Springer, 2005.

\bibitem{DF2} L.\ Dieulefait and N.\ Freitas
{\em Fermat-type equations of signature $(13,13,p)$ via Hilbert cuspforms},
Math. Annalen, 2013, (to appear).

\bibitem{ellenberg} J.\ Ellenberg
{\em Serre's conjecture over $\F_9$},
Ann. of Math. (2) {\bf 161} (2005), 1111--1142.


\bibitem{Elkies} N.\ Elkies,
{\em The Klein quartic in number theory},
in
{\em The Eightfold Way: The Beauty of Klein's Quartic Curve},
ed.\ Silvio Levy,
MSRI Publications
{\bf Volume 35}, 1998.

\bibitem{Faltings1} G.\ Faltings,
{\em Endlichkeitss\"{a}tze f\"{u}r abelsche Variet\"{a}ten \"{u}ber 
Zahlk\"{o}rpern},
Invent.\ Math.\ {\bf 73} (1983), no.\ 3, 349--366. 

\bibitem{Faltings} G.\ Faltings, 
{\em The general case of S.\ Lang’s conjecture}, 
pages 175–182 of Barsotti Symposium in Algebraic Geometry (Abano Terme, 1991), 
Perspect.\ Math.\ {\bf 15}, Academic Press, San Diego, CA, 1994.

\bibitem{recipes} N.\ Freitas,
{\em Recipes for Fermat-type equations of the form $x^r + y^r = Cz^p$},


\bibitem{Fisher} T.\ Fisher,
{\em The invariants of a genus one curves},
Proc.\ London Math.\ Soc.\ {\bf 97} (2008), 753--782.

\bibitem{FlynnTesta}
E.\ V.\ Flynn and D.\ Testa,
{\em Finite Weil restrictions of curves}, 
{\tt arXiv:1210.4407v3}, 19 May 2013.

\bibitem{FS} N.\ Freitas and S.\ Siksek,
{\em Modularity and the Fermat equation over totally real fields},
{\tt arXiv:1307.3162}, 11 July 2013.

\bibitem{Galbraith} S.\ D.\ Galbraith,
{\em Equations for Modular Curves}, DPhil thesis,
University of Oxford, 1996.

\bibitem{Gee} T.\ Gee,
{\em Automorphic lifts of prescribed types},
Mathematische Annalen {\bf 350} (2011), 107--144.


\bibitem{HK} 
E.\ Halberstadt and A.\ Kraus, 
{\em Sur la courbe modulaire $X_E(7)$}, 
Experiment.\ Math.\ {\bf 12} (2003), no.\ 1, 27--40.

\bibitem{HS} J.\ Harris and J.\ H.\ Silverman, 
{\em Bielliptic curves and symmetric products}, 
Proceedings of the American Mathematical Society {\bf 112} (1991), no.\ 2, 347--356.

\bibitem{Hartshorne} R.\ Hartshorne,
{\em Algebraic Geometry},
GTM {\bf 52}, Springer Verlag, 1977.


\bibitem{JM} F.\ Jarvis and J.\ Manoharmayum,
{\em On the modularity of supersingular elliptic curves
     over certain totally real number fields
},
Journal of Number Theory {\bf 128} (2008), no.\ 3, 589--618. 

\bibitem{JMee} F.\ Jarvis and P.\ Meekin,
{\em The Fermat equation over $\Q(\sqrt{2})$},
 Journal of Number Theory {\bf 109} (2004), no.\ 1, 182--196. 

\bibitem{Kamienny}
S.\ Kamienny, 
{\em Torsion points on elliptic curves and $q$-coefficients of modular forms},
Invent.\ Math.\ {\bf 109} (1992), no.\ 2, 221--229. 

\bibitem{KW1}
C.\ Khare and J.-P.\ Wintenberger, 
{\em Serre's modularity conjecture. I},
 Invent.\ Math.\ {\bf 178} (2009), no.\ 3, 485--504.

\bibitem{KW2}
C.\ Khare and J.-P.\ Wintenberger, 
{\em Serre's modularity conjecture. II},
 Invent.\ Math.\ {\bf 178} (2009), no. 3, 
505--586. 

\bibitem{KisM} H.\ Kisilevsky and M.\ R.\ Murty,
{\em Elliptic Curves and Related Topics},
CRM Proceedings and Lecture Notes {\bf 4}, 
American Mathematical Society, 1994.

\bibitem{Kisin} M.\ Kisin,
{\em Moduli of finite flat group schemes, and modularity},
 Annals of Mathematics. Second Series {\bf 170} (2009), no. 3,
  1085--1180.

\bibitem{antIII}
W.\ Kuyk and J.-P.\ Serre (editors),
{\em Modular Functions of One Variable III},
Lecture Notes in Mathematics {\bf 350},
Springer-Verlag, 1973.

\bibitem{Langlands} R.\ Langlands,
{\em Base Change for $GL(2)$},
Ann.\ of Math. Studies {\bf 96},
Princeton University Press, 1980.

\bibitem{Bao} B.\ Le Hung,
{\em Modularity of elliptic curves over some totally real fields},
{\tt arXiv:1309.4134}, 16 September 2013.

\bibitem{Ligozat} G.\ Ligozat,
{\em Courbes modulaires de niveau $11$},
pages 147--237 of \cite{antwerpV}. 

\bibitem{Mano1} J.\ Manoharmayum,
{\em On the modularity of certain $\GL_2(\F_7)$ Galois representations},
 Math.\ Res.\ Lett.\ {\bf 8} (2001), no.\ 5-6, 703--712. 

\bibitem{Mano2} J.\ Manoharmayum,
{\em Serre's conjecture for mod 7 Galois representations}
pages 141--149 of \cite{modcurves}.

\bibitem{Mazur} B.\ Mazur,
{\em Rational Points on Modular Curves},
pages 107--147 of \cite{antwerpV}.

\bibitem{MazurEisen}
B.\ Mazur,
{\em Modular curves and the Eisenstein ideal},
Inst.\ Hautes \'{E}tudes Sci.\ Publ.\ Math.\ {\bf 47} (1977), 33--186 (1978). 

\bibitem{Merel} L.\ Merel,
{\em Bornes pour la torsion des courbes elliptiques sur les corps de nombres},
Invent.\ Math.\ {\bf 124} (1996), no.\ 1-3, 437--449. 

\bibitem{Milne} J.\ S.\ Milne,
{\em Abelian Varieties}, Chapter V of \cite{CS}.


\bibitem{Muic} G.\ Mui\'{c},
{\em On embeddings of modular curves in projective space},
preprint, 17 March 2012.

\bibitem{PSS} B.\ Poonen, E.\ F.\ Schaefer and M.\ Stoll,
{\em Twists of $X(7)$ and primitive solutions to 
$x^2+y^3=z^7$},
Duke Math.\ J.\ {\bf 137} (2007), 103--158.

\bibitem{ogg} A.\ Ogg,
{\em Hyperelliptic modular curves},
 Bulletin de la Soci\'{e}t\'{e} Math\'{e}matique de France 
{\bf 102} (1974), 449--462.



\bibitem{Ribet} K.\ A.\ Ribet,
{\em Abelian varieties over $\Q$ and modular forms}
pages 241--261 of \cite{modcurves}.

\bibitem{Rohrlich2} D.\ E.\ Rohrlich,
{\em Elliptic curves and the Weil--Deligne group},
pages 125--158 of \cite{KisM}.

\bibitem{Rohrlich} D.\ E.\ Rohrlich,
{\em Modular curves, Hecke correspondences, and $\mathrm{L}$-functions},
pages 41--100 of \cite{Cornell}.

\bibitem{Rubin} K.\ Rubin,
{\em Modularity of mod $5$ representations},
pages 463--474.


\bibitem{Serre}
J.-P.\ Serre,
{\em Properi\'{e}t\'{e}s galoisiennes des points d'ordre fini
des courbes elliptiques},
Inventiones Math.\ {\bf 15} (1972), 259--331.



\bibitem{SeMordell}
J.-P.\ Serre,
{\em Lectures on the Mordell-Weil Theorem},
Aspects of Mathematics {\bf 15},
Vieweg, 
3rd edition, 1997.

\bibitem{SeGalois} J.-P.\ Serre,
{\em Topics in Galois Theory},
 A K Peters/CRC Press, 2nd edition, 2007.

\bibitem{antwerpV}
J.-P.\ Serre and D.\ B.\ Zagier (editors),
{\em Modular Functions of One Variable V},
Lecture Notes in Mathematics {\bf 601},
Springer-Verlag, 1977.


\bibitem{Siksek} S.\ Siksek,
{\em The modular approach to Diophantine equations},
pages 495--527 of \cite{Cohen}.

\bibitem{chabsym} S.\ Siksek,
{\em Chabauty for symmetric powers of curves},
Algebra \& Number Theory {\bf 3} (2009), no.\ 2, 209--236.

\bibitem{chabnf} S.\ Siksek,
{\em Explicit Chabauty over number fields}, 
Algebra \& Number Theory {\bf 7} (2013), no.\ 4, 765--793.

\bibitem{SilvermanII} J.\ H.\ Silverman,
{\em Advanced Topics in the Arithmetic of Elliptic Curves},
GTM 151, Springer, 1994.

\bibitem{SWI} C.\ Skinner and A.\ Wiles
{\em Residually reducible representations and modular forms},
Publications math\'ematiques de l'I.H.\'{E}.S.\ {\bf 89} (1999),
5--126.


\bibitem{SWII} C.\ Skinner and A.\ Wiles
{\em Nearly ordinary deformations of irreducible residual representations},
Ann.\ Fac.\ Sci.\ Toulouse Math.\ {\bf 10}, (2001), 185--215.



\bibitem{Stein} W.\ A.\ Stein,
{\em Modular Forms: A Computational Approach},
American Mathematical Society, 2007.

\bibitem{Stoll} M.\ Stoll,
{\em Implementing 2-descent for Jacobians of hyperelliptic curves},
Acta Arith.\ {\bf 98} (2001), 245--277.

\bibitem{StollRat} M.\ Stoll,
{\em Rational points on curves},
Journal de Th\'{e}orie des Nombres de Bordeaux {\bf 23} (2011), 257--277.

\bibitem{Sutherland} A.\ V.\ Sutherland,
{\em A local-global principle for
rational isogenies of prime degree},
Journal de th\'{e}orie des nombres de Bordeaux {\bf 24} (2012),
no.\ 2, 475--485.

\bibitem{SwD} H.\ P.\ F.\ Swinnerton-Dyer,
{\em On $\ell$-adic representations and congruences
for coefficients of modular forms},
pages 1--55 of \cite{antIII}.


\bibitem{TaylorHilbert}
R.\ Taylor, {\em On Galois representations associated to Hilbert
  modular forms}, Inventiones Mathematicae {\bf 98} (1989), no.~2,
  265--280.

\bibitem{TaylorIco} R.\ Taylor
{\em On icosahedral Artin representations II},
American Journal of Mathematics {\bf 125} (2003), 549--566.

\bibitem{TW} R.\ Taylor and A.\ Wiles,
{\em Ring-theoretic properties of certain Hecke algebras},
Annals of Mathematics {\bf 141} (1995), no.\ 3, 553--572.

\bibitem{Tunnell}
J.\ Tunnell,
{\em Artin's conjecture for representations of octahedral
type},
Bull.\ A.M.S.\ {\bf 5} (1981), 173--175.

\bibitem{Wiles} A.\ Wiles,
{\em Modular elliptic curves and Fermat's Last Theorem},
Annals of Mathematics {\bf 141} (1995), no.\ 3, 443--551.


\bibitem{Wiles2}
A.\ Wiles,
{\em On ordinary {$\lambda$}-adic representations associated to
modular forms}, Invent. Math.  \textbf{94} (1988), no.\ 3, 529--573.

\end{thebibliography}
\end{document}